\documentclass[12pt]{article}

\usepackage{amssymb,amsmath,xcolor,graphicx}
\begin{document}

%\begin{frontmatter}



\title{Periodic switching strategies for an isoperimetric control problem with application to nonlinear chemical reactions}


%\ead{benner@mpi-magdeburg.mpg.de}

%\ead{seidel@mpi-magdeburg.mpg.de}

\author{P.~Benner$^1$, A.~Seidel-Morgenstern$^1$, and A.~Zuyev$^{1,2}$\footnote{Corresponding author.}}
\date{\small $^1$Max Planck Institute for Dynamics of Complex Technical Systems, \\ Magdeburg, Germany \\
$^2$Institute of Applied Mathematics and Mechanics,\\ National Academy of Sciences of Ukraine}
%\cortext[cor1]{Corresponding author.}
%\ead{zuyev@mpi-magdeburg.mpg.de}
\maketitle

%\address{$^1$Max Planck Institute for Dynamics of Complex Technical Systems, \\ Sandtorstra{\ss}e 1, 39106 Magdeburg, Germany}

%\address{$^2$Institute of Applied Mathematics \& Mechanics,  \\ National Academy of Sciences of Ukraine}

\begin{abstract}
This paper deals with an isoperimetric optimal control problem for nonlinear control-affine systems with periodic boundary conditions.
As it was shown previously, the candidates for optimal controls for this problem can be obtained within the class of bang-bang input functions.
We consider a parametrization of these inputs in terms of switching times.
The control-affine system under consideration is transformed into a driftless system by assuming that the controls possess properties of a partition of unity.
Then the problem of constructing periodic trajectories is studied analytically by applying the Fliess series expansion over a small time horizon.
We propose analytical results concerning the relation between the boundary conditions and switching parameters for an arbitrary number of switchings.
These analytical results are applied to a mathematical model of non-isothermal chemical reactions.
It is shown that the proposed control strategies can be exploited to improve the reaction performance in comparison to the steady-state operation mode.
\end{abstract}

%\begin{keyword}
%isoperimetric optimal control problem \sep Fliess series \sep bang-bang control \sep chemical reaction


%\MSC[2010] 93C10 \sep 	93C95 \sep 49K15 \sep 	49N20


%\end{keyword}

%
%\end{frontmatter}

\newpage
\noindent{\bf Highlights}
\begin{itemize}
\item Bang-bang control strategies corresponding to periodic trajectories of nonlinear control-affine systems with isoperimetric constraints are studied in detail.
\item An approximate method for defining the switching parameters in terms of solutions to auxiliary algebraic equations is proposed.
\item The proposed control strategies are applied to the optimal control of non-isothermal chemical reactors.
\item The performance improvement of a nonlinear reaction is confirmed by analytical and numerical results.
\end{itemize}


\section{Introduction}\label{sec_1}

The issues of optimal design of processes in chemical engineering give rise to a series of challenging problems in mathematical control theory.
An important class of these problems is related to the analysis of mathematical models of nonlinear chemical reactions described by ordinary differential equations under periodic control strategies (see, e.g,~\cite{SH2013,NSP2015,NP2016} and references therein).
It was pointed out in~\cite{Douglas1967} that the performance measure of a nonlinear chemical reaction essentially depends on characteristics of the periodic input signal corresponding to the feed composition.
{
A particular estimate of the improvement of the objective function for different types of periodic controls was proposed in~\cite{SB1980} by computing the second variation of the cost in the frequency domain.
It was shown there that small harmonic oscillations in feed reactant concentration lead to the increase of the mean  product concentration for a class of second-order isothermal chemical reactions.
}
Over the last few decades, the optimization problems for mathematical models of periodic chemical reactions have received a lot of attention in both mathematical and engineering studies.

{
Sufficient conditions for improving the performance of a periodic process in a neighborhood of the given steady-state were proposed in~\cite{BH1971}.
These conditions are based on the Pontryagin maximum principle and relaxed steady-state analysis.
It was noted that the relaxed steady-state analysis can indicate that the improvement of large magnitude is possible by cycling, while the maximum principle provides no information about the magnitude of possible improvement.
}

{
The Laplace--Borel transform has been applied in~\cite{HP1997} for the study of forced nonlinear processes with hyperbolic equilibrium points.
This approach is equivalent to the transfer function approach with a modified set of zeros and poles to mimic nonlinearities in the system. The modified transfer function is computed in~\cite{HP1997} for a second-order isothermal reaction to obtain minimum and maximum bounds for the amplitude and phase to single-tone inputs.
}
{
Another method for finding periodic trajectories corresponding to the maximum of the time-average output of a nonlinear system  was proposed in~\cite{GDPH2007} within the framework of extremum seeking control.
By assuming that the considered control system is flat, the original optimal control problem is transformed to a parameterized optimization task.
This idea is applied to a drug delivery example in~\cite{GDPH2007}.
}

{
Second-order necessary optimality conditions have been studied in~\cite[Chapter~IX]{C1988} for periodic problems described by nonlinear ordinary differential equations with state constraints and an isoperimetric condition.
A $\Pi$-test under state constraints is proposed to estimate optimal periodic solutions in a neighborhood of the steady-state of the system.
This result extends the second variation technique in the frequency domain developed in~\cite{SB1980} for unconstrained problems.
As an example, the $\Pi$-test is applied to a non-isothermal chemical reaction with a single control in~\cite{C1988}.
The considered control is proportional to the overall heat transfer coefficient and can be adjusted by the coolant flow rate.
The behavior of this dynamical system in a neighborhood of its steady-state is analyzed by treating the Damk{\"oh}ler number as a bifurcation parameter.
}

{
A generalized $\Pi$-criterion was proposed in the paper~\cite{SY1990} to describe the effects of periodic perturbations around any given steady-state of a nonlinear system. This criterion is not restricted just to the optimal steady-states, but allows studying an arbitrary equilibrium point under the stability assumption. The criterion is applied to a nonlinear system describing two parallel chemical reactions of the type $2A \to B$ and $A\to C$.
This methodology has been further developed in~\cite{SY1991} for other types of chemical reactions, in particular, $A\to B\to C$. The generalized $\Pi$-criterion has been applied to find constraints on the activation energy such that the yield of the product is improved using high-frequency periodic perturbations in the temperature.
The $\Pi$-criterion was also successfully applied to broad classes of nonlinear systems with multidimensional inputs, including mathematical models of continuous bioprocesses governed by the conservation equations for cell mass. A particular class of such models governed by third-order nonlinear ordinary differential equations with integral constraints was considered in~\cite{P2000}.
The differential equations of this type have been further addressed in the paper~\cite{ZSP2017} by exploiting the Laplace--Borel transform and separating the stationary behavior and the transient dynamics.
The above study resulted in the computation of the proper forcing amplitude and frequency to optimize the performance measure.
}

{
The $\Pi$-criterion was adapted to periodic inputs with very low frequencies and tested on an isothermal continuous stirred-tank reactor (CSTR) in~\cite{P2003}.
For a set of non-isothermal parallel reactions in a CSTR, the $\Pi$-criterion was applied in~\cite{LM2001}.
In the cited paper, a numerical simulation is performed with periodic forcing of the feed temperature, and the optimal frequency for the best performance improvement is found numerically.
It is noted in~\cite{KDT2012} that the $\Pi$-criterion has local nature and provides approximate estimates of the performance measure for the inputs with small amplitudes.
Thus, higher-order corrections to the $\Pi$-criterion are proposed in~\cite{KDT2012}, based on approximations of the center manifold by power series.
As a result, a truncated series is obtained to compute approximately the performance measure of a nonlinear control problem with harmonic inputs.
This approximation is applied to a non-isothermal CSTR controlled by the temperature variation under sinusoidal control strategies.
The case of square wave inputs is considered by using the truncated Fourier series. It is shown that the use of higher-order terms can improve the accuracy of analytic estimates of the performance under large control amplitudes.
}

{
A method for solving an optimal control problem within the class of periodic inputs based on the Carleman linearization was proposed in~\cite{HLPS1993}.
The authors of the cited paper considered pulsed periodic modulations of single-input nonlinear systems and applied their construction to optimize an isothermal CSTR.
It was shown that the second-order Carleman procedure produces a good approximation of the results from direct numerical integration.
}

{
Forced oscillations in an exothermic CSTR were studied in~\cite{CFMS1987} within the framework of vibration control.
This approach resulted in a modification of the dynamical properties of the nonlinear system under consideration by using fast oscillations of the input flow rate.
Such input modulations are designed to ensure an asymptotically stable periodic operation of the reactor in a neighborhood of its unstable equilibrium.
Analytic stability conditions have been derived by using the averaging method, and numerical simulations have been performed together with experiments to verify the stability of the controlled system.
A rather general version of the averaging method for multifrequency systems with exact estimates of solutions in powers of a small parameter was presented in~\cite{KMBMPSZS2009}.
The proposed method also allows studying conditional asymptotic stability of integral manifolds and averaging the boundary conditions together with the differential equations.
Note that the averaging technique was previously applied to nonlinear systems with Arrhenius-type dynamics in the paper~\cite{BBM1983}.
It was shown there that an increase in productivity of catalytic reactors, modeled as Arrhenius systems, can be achieved by means of vibrational stabilization.
}

{
Recent results on estimating the time-average performance of chemical reactions by the nonlinear frequency response method with periodic inputs can be found in~\cite{PNS2018,SLZYW2018}.
}

In the paper~\cite{CES2017}, the problem of maximizing the performance of a periodic chemical reactor is treated as an optimal control problem with isoperimetric constraints.
It has been shown that each optimal control for this problem is a bang-bang control,
and the maximal number of switchings has been estimated for the linearized equations.
Although some basic properties of the extremal controls have been analysed in~\cite{CES2017,DSTA2017} for a particular first-order non-isothermal reaction, the general question of computing control strategies for problems with isoperimetric constraints and periodic boundary conditions remains open.
{
The reported few attempts to apply experimentally the principle of forced periodic operation of chemical reactors were recently summarized in~\cite{SH2013}.
}

The objective of our work is to propose an efficient technique for defining the switching controls for a wide class of nonlinear control-affine systems with isoperimetric constraints, and to apply these theoretical results to mathematical models of nonlinear chemical reactions with periodic inputs.
The current paper originates from our previous conference paper~\cite{DSTA2017} and essentially extends the analytical approach of~\cite{DSTA2017} for the case of bang-bang controls with an arbitrary number of switchings.
The main theoretical results are presented in Section~3 and tested by numerical simulations in Section~4.
We present a novel class of conditions linking together the switching parameters and the initial data. Such conditions are formulated as systems of algebraic equations involving vector fields of the system and their Lie derivatives at the initial point.
The efficiency of the control design scheme, based on these algebraic equations,
is illustrated with an example of a non-isothermal chemical reaction of the type ``$A\to$ product''.

\section{Mathematical model of a non-isothermal chemical reaction}
\label{section2}
Our study is motivated by optimal control problems for mathematical models of non-isothermal chemical reactions governed by nonlinear ordinary differential equations.
As an important representative of this class of models, we consider a simple reaction of the type ``$A\to$ product''
described by the following control-affine system (see~\cite{CES2017,DSTA2017}):
\begin{equation}
\dot x = f_0(x) + u_1 g_1(x)+ u_2 g_2(x),\quad x =\begin{pmatrix}x_1 \\ x_2\end{pmatrix} \in{\mathbb R}^2,\; u = \begin{pmatrix}u_1 \\ u_2\end{pmatrix} \in U \subset {\mathbb R}^2,
\label{csys}
\end{equation}
$$
U= [u_1^{min},u_1^{max}]\times [u_2^{min},u_2^{max}].
$$
The components of the state vector $x(t)$ have the following physical meaning: $x_1(t)$ describes to the outlet concentration of $A$,  and $x_2(t)$ corresponds to the temperature in the reactor at time $t$.
System~\eqref{csys} is controlled by modulating the inlet concentration of $A$ (input signal $u_1(t)$)
and the temperature of the inlet stream (input signal $u_2(t)$).
The values of $x_1$, $x_2$, $u_1$, and $u_2$ are taken as dimensionless deviations of the corresponding physical quantities from their steady-state values under a suitable rescaling.
Thus, system~\eqref{csys} admits the equilibrium $x_1=x_2=0$ with $u_1=u_2=0$ which corresponds to a certain operating mode of the reactor with constant inflow and outflow characteristics.

For the mathematical model of the $\bar n$-th order reactions considered in~\cite{CES2017}, the vector fields of system~\eqref{csys} are
\begin{equation}
f_0(x) = \begin{pmatrix}k_1 e^{-\varkappa}-\phi_1 x_1 - k_1 (x_1+1)^{\bar n} e^{-\varkappa/(x_2+1)} \\ k_2 e^{-\varkappa}-\phi_2 x_2 - k_2 (x_1+1)^{\bar n} e^{-\varkappa/(x_2+1)}\end{pmatrix}, \; g_1 = \begin{pmatrix}1 \\ 0 \end{pmatrix},\;
g_2 = \begin{pmatrix}0 \\ 1 \end{pmatrix},
\label{f_notations}
\end{equation}
where $\varkappa$, $k_i$, and $\phi_i$ are parameters of the reaction.
All the necessary details concerning the derivation of control system~\eqref{csys} can be found in~\cite{NSP2015,NP2016,CES2017}.

As it has been already noted in the introduction, a series of challenging problems is related to the optimization of periodic operating modes for nonlinear chemical reactions. One of this problems deals with maximizing the rate of conversion of $A$ to the final product.
 The equivalent task of minimizing the mean concentration of $A$ at the output of the reactor has been considered in~\cite{CES2017,DSTA2017} in the context of isoperimetric optimal control problems.

%To formulate this problem, we introduce the class of admissible controls ${\cal U}_{{\tau}}$ consisting of all measurable functions $u:[0,{\tau}]\to U\subset {\mathbb R}^2$.
\subsection{Isoperimetric optimal control problem}
\label{subsection21}

We recall the problem statement from~\cite{CES2017} below.

{\bf Problem~2.1.}
{\em
For given ${\tau}>0$, $x^0\in \mathbb R^2$, and $\bar u_1\in \mathbb R$, the goal is to find a control $\hat u\in L^\infty([0,{\tau}];U)$ that minimizes the cost
$$
J[x]= \frac{1}{{\tau}}\int_0^{{\tau}} x_1(t)dt
$$
along the solutions  $x(t)$ of system~\eqref{csys} corresponding to the admissible controls $u\in L^\infty([0,{\tau}];U)$
such that the periodic boundary conditions
$$
x^0=x(0)=x({\tau})
$$
and the isoperimetric constraint
$$
\frac{1}{{\tau}}\int_0^{{\tau}} u_1(t)dt = \bar u_1
$$
hold.
}

This problem statement has clear physical meaning: to minimize the remainder of $A$ in the output of the reactor by consuming a fixed amount of $A$ over the period of time $t\in [0,{\tau}]$.
As one can easily see, the cost $J$ vanishes on the steady-state solution $x(t)\equiv $ with $u(t)\equiv 0$.
So, in particular, if one can construct a ${\tau}$-periodic control $u(t)$ such that $\int_0^{{\tau}} u_1(t)dt=0$ and the corresponding solution $x(t)$ is ${\tau}$-periodic with $
J[x]<0$, then such control improves the performance of the reaction with respect to its steady-state operation.

{\bf Remark~2.1.}
{\em
The above isoperimetric problem was analyzed in~\cite{CES2017} by using a modification of the Pontryagin maximum principle with Lagrange multipliers (cf.~\cite{Schmitendorf,C1988}).
It was shown in~\cite{CES2017} that, if $\hat u(t)$ is an optimal control for the above problem,
then the minimizer $\hat u(t)$ can be chosen in the class of bang-bang controls.
The number of switchings of $\hat u(t)$ was also estimated in~\cite{CES2017}
for the linearization of system~\eqref{csys} with ${\bar n}=1$, when
the differential equations for $x(t)$ and the
adjoint variables $p(t)$ of the Hamiltonian system are decoupled.
In the latter case, it was shown that any bang-bang control $\hat u(t)$ satisfying the Pontryagin maximum principle has at maximum 4 switchings in the interval $t\in[0,{\tau}]$, provided that
\begin{equation}
 D=(\phi_1+\phi_2+{\bar n} \tilde k_1 + \varkappa \tilde k_2 )^2 - 4(\phi_1 \phi_2 + \phi_1 \varkappa \tilde k_2 + {\bar n} \phi_2 \tilde k_1)>0,
 \label{D_cond}
\end{equation}
where $\tilde k_1 =  k_1 e^{-\varkappa}$, $\tilde k_2 =  k_2 e^{-\varkappa}$.
Note that the condition~\eqref{D_cond} holds for the reactor model considered in~\cite{CES2017}.}

Based on the assumption $D>0$, some particular switching strategies for Problem~2.1 have been analysed in the previous papers~\cite{CES2017,DSTA2017} for the the control system~\eqref{csys} with the vector fields given by~\eqref{f_notations}.
We will consider this problem in a more general form below.

\section{Main results}
\label{section3}
In this section, we will address the problem of defining the bang-bang controls with desired properties
 for general control-affine systems of the form
 \begin{equation}
\dot x = f_0(x) + \sum_{j=1}^m u_j g_j(x),\quad x=(x_1,...,x_n)^T\in X,\; u=(u_1,...,u_m)^T \in U,
\label{controlaffine}
\end{equation}
where the vector fields $f_0$, $g_1$, ..., $g_m$ are assumed to be smooth in the domain $X\subset \mathbb R^n$, and the set of control values $U \subset \mathbb R^m$ is compact.
We also assume that $(0,0)\in X\times U$ and $f_0(0)=0$, so that system~\eqref{controlaffine} admits the equilibrium $x=0$ with $u=0$.
For a time horizon ${\tau}>0$, we denote the set of admissible controls as ${\cal U}_{{\tau}}=\{u\in L^\infty[0,{\tau}]: u(t)\in U\; \text{for all}\;t\in [0,{\tau}]\}$.
In the sequel, we will treat the control system~\eqref{csys} as a particular case of system~\eqref{controlaffine} and study properties of related switching control strategies for arbitrary dimensions $n$ and $m$.

\subsection{Bang-bang inputs for general control-affine systems}

%for the corresponding solution $x(t)$ of system~\eqref{controlaffine}.
% Starting from a particular two-dimensional system considered in~\cite{DSTA2017}, we will analyse the
%For the rest of this paper, we assume that the number of switchings is finite for each optimal control $\tilde u(t)$.

In order to study the family of bang-bang controls satisfying the isoperimetric constraint
\begin{equation}
\frac{1}{{\tau}}\int_0^{{\tau}} u_1(t)dt = \bar u_1
\label{constraint}
\end{equation}
such that the corresponding solution $x(t)$ of system~\eqref{controlaffine} satisfies
\begin{equation}
x^0=x(0)=x({\tau}),
\label{periodic_BC}
\end{equation}
we fix an integer $N\ge 1$ and consider a finite sequence of control values
\begin{equation}
u^1,u^2,...,u^N \in \partial U
\label{strategyG}
\end{equation}
together with a partition of the time interval
\begin{equation}
0=t_0 < t_1< ... < t_{N} = {\tau},\quad  (\tau_j=t_j - t_{j-1}>0),
\label{partition}
\end{equation}
in order to define the following bang-bang control $u\in {\cal U}_{{\tau}}$:
\begin{equation}
u(t)=u^j \quad \text{for}\;\; t\in [t_{j-1},t_{j}),\quad j=\overline{1,N}.
 \label{u_tilde}
\end{equation}
%The above procedure provides a parametrization of the bang-bang controls in terms of switching scenarios.

As it follows from the results of Section~\ref{section2} for a particular form of system~\eqref{controlaffine} with $n=2$ and $U= [u_1^{min},u_1^{max}]\times [u_2^{min},u_2^{max}]$, {\em if $u\in {\cal U}_{{\tau}}$ is an optimal control for Problem~2.1,
then $u(t)$ can be constructed in the form~\eqref{u_tilde}
by minimizing the corresponding cost
\begin{equation}
J[x]= \frac{1}{{\tau}}\int_0^{{\tau}} x_1(t)dt
\label{J_cost}
\end{equation}
for all possible choices of $N$, switching scenarios~\eqref{strategyG}, and switching times~\eqref{partition}, such that the constraints~\eqref{constraint}--\eqref{periodic_BC} are satisfied.
}

%For this purpose we first embed the isoperimetric constaint~\eqref{constraint} into the set of periodic boundary conditions by introducing an additional differential equation $\dot \xi_{n+1}=u_1-\bar u_1$.
%If $\xi_{n+1}(t)$ is a solution of this equation with a control $u(t)$ on $[0,{\tau}]$, then~\eqref{constraint} is equivalent to the periodic boundary condition $\xi_{n+1}(0)=\xi_{n+1}({\tau})$.
%Thus the problem of defining an admissible control $u\in{\cal U}_{{\tau}}$ such that~\eqref{constraint}--\eqref{periodic_BC} holds may be considered as the problem with periodic boundary conditions for the extended system in the $(n+1)$-dimensional space.

\subsection{Reduction to driftless systems and the Fliess expansion}
Our goal is to propose an efficient control design scheme that allows computing the switching parameters $\tau_1$, $\tau_2$, ..., $\tau_N$ from auxiliary algebraic equations.
For this purpose we first rewrite~\eqref{controlaffine} with the controls~\eqref{u_tilde} as the driftless system
\begin{equation}
\dot x = \sum_{j=1}^N v_j(t) f_j(x)
\label{driftless}
\end{equation}
with
\begin{equation}
f_j(x)=f_0(x)+\sum_{i=1}^m u^j_i g_i(x)
\label{f_driftless}
\end{equation}
and
\begin{equation}
v_j(t)=\chi_{[t_{j-1},t_j)}(t) \quad \text{for}\; j=1,2,...,N.
\label{vj}
\end{equation}
 Here $\chi_{[t_{j-1},t_j)}(t)$ is the indicator function: $\chi_{[t_{j-1},t_j)}(t)=1$ if $t\in [t_{j-1},t_j)$,
and  $\chi_{[t_{j-1},t_j)}(t)=0$ if $t\notin [t_{j-1},t_j)$.
The family of control functions $\{v_1,v_2,...,v_N\}$ possesses an important property of a partition of unity: $\sum_{j=1}^N v_j(t)=1$ for all $t\in [0,{\tau})$.

{
For each vector field $f_j$, $j=1,2,...,N$, we denote its flow on $X$ by $e^{tf_j}$,
i.e. $x(t)=e^{tf_j}(x^0)$ stands for the solution to the Cauchy problem $\dot x(t)=f_j(x(t))$ with the initial value $x(0)=x^0\in X$.
Then the problem of finding an admissible bang-bang control of the form~\eqref{u_tilde} such that the corresponding solution of~\eqref{controlaffine} satisfies~\eqref{periodic_BC} can be formulated as the problem of defining positive numbers $(\tau_1,\tau_2,...,\tau_N)$ such that
\begin{equation}
e^{\tau_N f_N} \circ ...\circ e^{\tau_2 f_2} \circ e^{\tau_1f_1}(x^0)=x^0
\label{periodic_exp}
\end{equation}
or, equivalently,
\begin{equation}
e^{\tau_i f_i} \circ ... \circ e^{\tau_1f_1}(x^0) = e^{-\tau_{i+1} f_{i+1}} \circ ...\circ e^{-\tau_N f_N}(x^0) \quad \text{for some}\;i< N.
\label{periodic_expm}
\end{equation}
From the geometric viewpoint, the study of admissible periodic trajectories for Problem~2.1 is thus reduced to the construction of closed curves in $X$ from arcs of the type $e^{tf_j}(x)$.
Although the two conditions~\eqref{periodic_exp} and~\eqref{periodic_expm} are equivalent, the application of~\eqref{periodic_expm} may have an advantage in computation if the $e^{\tau_j f_j}(x)$ are growing fast for large values of $\tau_j$.
In the sequel, we exploit formula~\eqref{periodic_exp} as a condition between the initial value $x^0$ and switching times $(\tau_1,\tau_2,...,\tau_N)$ if  the time horizon $\tau>0$ is small enough. As the left-hand side of~\eqref{periodic_exp} is the solution of~\eqref{driftless},~\eqref{vj} at time $\tau$, we will apply the Fliess functional expansion for formal manipulations with such solutions.
}

If the vector fields $f_j$ are analytic and $y=h(x)$ is an analytic output function,
then the value of $y(t)=h(x(t))$ for the corresponding solution $x(t)$ of system~\eqref{driftless}  with the initial data $x(0)=x^0$
admits the following representation~\cite{Lamna},~\cite[Chapter~4]{NS90}:
\begin{equation}
y(t) = h(x^0) + \sum_{\nu=0}^\infty \sum_{i_1,...,i_\nu=1}^N L_{f_{i_1}} \cdots L_{f_{i_\nu}} h(x^0) \int_0^t d\xi_{i_\nu} \cdots d\xi_{i_1}, \quad t\in [0,{\tau}],
 \label{Fliess}
\end{equation}
where $L_{f_i}h(x) = \frac{\partial h(x)}{\partial x}f_i(x)$ denotes the Lie derivative, and $\frac{\partial h(x)}{\partial x}$ is the Jacobian matrix.
{Note that the Fliess expansion~\eqref{Fliess} can be obtained from the Volterra series (see, e.g.,~\cite{Lamna}).}
The iterated integrals in~\eqref{Fliess} are defined as $\int_0^t d\xi_i = \xi_i(t) = \int_0^t v_i(t)dt$ for $i=1,2,...,N$, and, by induction,
$$
\int_0^t d\xi_{i_\nu} \cdots d\xi_{i_1} = \int_0^t d\xi_{i_\nu}(s) \int_0^s d\xi_{i_{\nu-1}} \cdots d\xi_{i_1}.
$$
% The Fliess series~\eqref{Fliess} converge uniformly for small $t$ and small $|u_i(t)|$.

In particular, the first terms of the expansion~\eqref{Fliess} for the controls given by~\eqref{vj} and $h(x)=x$ can be written as
\begin{equation}
\scriptsize
x(t) = x^0 + \sum_{i=1}^N   f_i(x^0) V_i(t) + \sum_{i,j=1}^N\bigl(L_{f_j} f_i \bigr)(x^0) V_{ij}(t)  + \sum_{i,j,l=1}^N  \bigl( L_{f_l} L_{f_j} f_i \bigr) (x^0)V_{ijl}(t) + R(t),\;t\in [0,{\tau}],
\label{Fliess3}
\end{equation}
where
\begin{equation}
\begin{aligned}
V_i(t)&=\int_0^t v_i(s )ds,\\ V_{ij}(t)&=\int_0^t \int_0^s v_i(s)v_j(p)dp\,ds,\\ V_{ijl}(t)&=\int_0^t \int_0^s \int_0^p v_i(s)v_j(p)v_l(r)dr\, dp\,ds.
\end{aligned}
\label{Vall}
\end{equation}
Note that the expansion~\eqref{Fliess3} is valid if the vector fields $f_j$ are of class $C^3(X)$, and its remainder admits the estimate $\max_{[0,{\tau}]}\|R(t)\|=O({\tau}^4)$, see~\cite{ECC2016,IJC2017} {for the proof}.
Throughout this paper we use the asymptotic notation $O(t^k)$ for small values of $t>0$: we write  $\phi(t)=O(t^k)$
 if and only if $\limsup_{t \downarrow 0} |\phi(t)|t^{-k}< \infty$.

\subsection{Control design scheme}

%As it was noted in Section~3.1, the problem of defining the switching times

The basic result we will prove is as follows.

{\bf Theorem~3.1.}
{\em
Let $u(t)$ be a control defined by~\eqref{u_tilde} with some parameters $N\ge 1$, $\{u^1,u^2,...,u^N\}\subset U$, $0=t_0<t_1<t_2<...<t_N={\tau}$,
and let $x(t)$, $t\in [0,{\tau}]$, be the corresponding solution of system~\eqref{controlaffine} with an initial data $x(0)=x^0\in X$.
If the conditions~\eqref{constraint} and~\eqref{periodic_BC} are satisfied, then
\begin{equation}
\sum_{i=1}^N \tau_i u^i_1 = {\tau} \bar u_1,\; \tau_i = t_i-t_{i-1}>0,
\label{isotau}
\end{equation}
and
\begin{equation}
\sum_{i=1}^N \left( \tau_i f_i +\frac{\tau_i^2}{2}  L_{f_i}f_i \right) + \sum_{1\le j<i\le N} \tau_i \tau_j L_{f_j}f_i + \sum_{i,j,l=1}^N \int_{t_{l-1}}^{t_l}V_{ij}(t)dt\, L_{f_j} L_{f_i}f_l  = O({\tau}^4),
\label{periodic_tN}
\end{equation}
where
\begin{equation}
2V_{ii}(t)=\left\{\begin{array}{ll} 0,&t\le t_{i-1},\\ {(t-t_{i-1})^2},&t\in (t_{i-1},t_i), \\ \tau_i^2, & t\ge t_i,\end{array}\right.\quad V_{ij}(t)=0\;\;\text{for}\; i<j,
\label{Vii}
\end{equation}
\begin{equation}
V_{ij}(t)=\left\{\begin{array}{ll} 0,&t\le t_{i-1},\\ (t-t_{i-1})\tau_j,&t\in (t_{i-1},t_i), \\ \tau_i \tau_j, & t\ge t_i,\end{array}\right. \quad \text{for} \; i>j.
\label{Vij}
\end{equation}
Moreover, the cost~\eqref{J_cost} is equal to $J[x]={\bar x}_1$, where
\begin{equation}
\footnotesize
\begin{aligned}
\bar x & = \frac{1}{{\tau}}\int_0^{{\tau}}x(t)dt = x^0 + \frac{1}{2{\tau}}\sum_{i=1}^N \tau_i\left(\tau_i+2({\tau}-t_i)\right)f_i \\
& + \frac{1}{6 {\tau}}\sum_{i=1}^N \tau_i^2 \left(\tau_i + 3 ({\tau}-t_i)\right) L_{f_i} f_i + \frac{1}{2{\tau}}\sum_{1\le j<i\le N} \tau_i \tau_j \left( \tau_i + 2 ({\tau}-t_i)\right) L_{f_j} f_i + O({\tau}^3).
\end{aligned}
\label{barx}
\end{equation}
The vector fields in formulas~\eqref{periodic_tN} and~\eqref{barx} are evaluated at $x=x^0$.
}

{\bf Proof.}
If the piecewise-constant control $u(t)$ is given by formula~\eqref{u_tilde} then $\int_0^{{\tau}}u(t) dt = \sum_{i=1}^N \tau_i u^i$, and the isoperimetric constraint~\eqref{constraint} is reduced to~\eqref{isotau}.
Let $x(t)$ be the solution of system~\eqref{controlaffine} corresponding to the initial data $x(0)=x^0$ and $u=u(t)$, then $x(t)$ is also a solution of system~\eqref{driftless} with the control~\eqref{vj}, so that we will use the Fliess expansion~\eqref{Fliess3} to prove formulas~\eqref{periodic_tN} and~\eqref{barx}.
Straightforward computation of the integrals in~\eqref{Vall} yields
\begin{equation}
V_i(t)=\left\{\begin{array}{ll} 0,&t\le t_{i-1},\\ t-t_{i-1},&t\in (t_{i-1},t_i), \\ \tau_i, & t\ge t_i,\end{array}\right.
\label{Vi}
\end{equation}
together with the relations~\eqref{Vii},~\eqref{Vij}, and
\begin{equation}
V_{ijl}(t)=\int_0^t v_i(s)V_{jl}(s)ds,\; V_{ijl}({\tau})= \int_{t_{i-1}}^{t_i} V_{jl}(s)ds.
\label{Vijl}
\end{equation}
Then the assertion~\eqref{periodic_tN} follows from the boundary condition $x(0)=x({\tau})$ and formulas~\eqref{Fliess3},~\eqref{Vii},~\eqref{Vij},~\eqref{Vi},~\eqref{Vijl} with $t={\tau}$.
Similarly we obtain the representation~\eqref{barx} by expressing $\int_0^{{\tau}}x(t)dt$ from~\eqref{Fliess3}.
$\square$

We also deduce one corollary of Theorem~3.1 with the following parametrization of switching times:
\begin{equation}
\begin{aligned}
\tau_j & = \alpha_j {\tau},\; \alpha_j>0,\quad j=2,..., N, \\
\tau_1 & = \left(1-\sum_{j=2}^N \alpha_j\right){\tau}>0.
 \end{aligned}
\label{tau-alpha}
\end{equation}
The above notations are convenient for eliminating $\tau_1$ from~\eqref{periodic_tN} and~\eqref{barx}.
Then Theorem~3.1 implies

{\bf Corollary~3.1.}
{\em
Let $u(t)$ be a control defined by~\eqref{u_tilde} with some parameters $N\ge 2$, $\{u^1,u^2,...,u^N\}\subset U$, $0=t_0<t_1<t_2<...<t_N={\tau}$,
and let $x(t)$, $t\in [0,{\tau}]$, be the corresponding solution of system~\eqref{controlaffine} with an initial data $x(0)=x^0\in X$.
If the conditions~\eqref{constraint} and~\eqref{periodic_BC} are satisfied, then
\begin{equation}
{
u_1^1+\sum_{i=2}^N \alpha_i (u^i_1 -u^1_1)
}= \bar u_1,
\label{isotau-alpha}
\end{equation}
and
\begin{equation}
\footnotesize
\begin{aligned}
& f_1+\sum_{j=2}^N \alpha_j (f_j-f_1) + \frac{{\tau}}{2}\Bigl\{L_{f_1}f_1 + 2\sum_{j=2}^N \alpha_j L_{f_1}(f_j-f_1)  \Bigr.\\
& \Bigl. + \sum_{j=2}^N \alpha_j^2(L_{f_1}(f_1-2f_j) + L_{f_j}f_j) + 2\sum_{2\le j<i\le N} \alpha_i \alpha_j (L_{f_1}(f_1-f_i-f_j)+ L_{f_j}f_i)\Bigr\} = O({\tau}^2),
\end{aligned}
\label{periodicN5}
\end{equation}
\begin{equation}
\begin{aligned}
\bar x = & x^0 + \frac{{\tau}}{2}\left\{f_1 + \sum_{j=2}^N \alpha_j \left(\alpha_j + 2 \sum_{i=j+1}^N\alpha_i\right)(f_j-f_1)\right\}\\
& +\frac{{\tau}^2}{6}\Bigl\{L_{f_1}f_1+ 3\sum_{j=2}^N \alpha_j\left(\alpha_j+ 2 \sum_{i=j+1}^N\alpha_i \right) L_{f_1}(f_j-f_1)\Bigr. \\
& \Bigl. + \sum_{j=2}^N \alpha_j^2\left(\alpha_j + 3\sum_{i=j+1}^N\alpha_i\right) (L_{f_1}(2 f_1 -3 f_j) + L_{f_j}f_j) \Bigr. \\
& \Bigl. + 6\sum_{2\le i<j<l\le N}\alpha_i\alpha_j\alpha_l (L_{f_1}(2f_1 - 2f_i - f_j)+ L_{f_i}f_j)\Bigr. \\
& \Bigl. + 3\sum_{2\le i<j\le N}\alpha_i\alpha_j^2(L_{f_1}(2f_1 - 2f_i-f_j)+L_{f_i}f_j)\Bigr\} + O({\tau}^3),
\end{aligned}
\label{barxN5}
\end{equation}
where $\alpha_i>0$ are related to $\tau_i=t_i-t_{i-1}>0$ by means of~\eqref{tau-alpha}.
}

{
We also formulate particular corollaries of the above result for the case $\bar u_1=0$ and symmetric $U\subset {\mathbb R}^m$, i.e., if $u\in U$ implies $-u\in U$.
}

{
{\bf Corollary~3.2.}
{\em Let the assumptions of Corollary~3.1 be satisfied with $N=2$, $\bar u_1=0$, $u_1^1\neq 0$, and let $\tau_1,\tau_2>0$ be related to $\alpha_2$ by means of~\eqref{tau-alpha}.
Then
\begin{equation}
u^2_1 = - u^1_1,\;\tau_1=\tau_2=\frac{\tau}{2},\quad \alpha_2=\frac{1}{2},
\label{tauN2}
\end{equation}
and
\begin{equation}
f_1 + f_2 + \frac{{\tau}}{4}\left(L_{f_1}f_1 - L_{f_2}f_2 \right)+ \frac{{\tau^2}}{24}\left( L_{f_1}L_{f_1}f_1 + L_{f_2}L_{f_2}f_2\right)= O({\tau}^2),
\label{periodicN2}
\end{equation}
\begin{equation}
\bar x = x^0 + \frac{{\tau}}{8}\left(f_1 -f_2\right) + \frac{{\tau}^2}{48}\left(L_{f_1}f_1 + L_{f_2}f_2\right) + O({\tau}^3).
\label{barxN2}
\end{equation}
}
}

{
{\bf Corollary~3.3.}
{\em Let the assumptions of Corollary~3.1 be satisfied with $N=3$, $\bar u_1=0$, $u_1^1=-u^2_1=-u^3_1\neq 0$, and let $\tau_1,\tau_2,\tau_3>0$ be related to $\alpha_2, \alpha_3$ by means of~\eqref{tau-alpha}.
Then
\begin{equation}
\tau_1=\frac{{\tau}}{2},\; \tau_2 = \alpha_2 {\tau},\; \tau_3 = \left(\frac{1}{2}-\alpha_2\right)\tau,\quad \alpha_2\in \left(0,\frac{1}{2}\right),
\label{tauN3}
\end{equation}
and
\begin{equation}
\small
\begin{aligned}
f_1&+f_3 + 2\alpha_2(f_2-f_3)
+ \frac{{\tau}}{4}\Bigl\{  L_{f_1}f_1+2L_{f_1}f_3 +L_{f_3}f_3 \Bigr. \\
&\Bigl.+4\alpha_2(L_{f_1}f_2-L_{f_1}f_3+L_{f_2}f_3-L_{f_3}f_3) +4\alpha_2^2 (L_{f_2}f_2-2L_{f_2}f_3+L_{f_3}f_3) \Bigr\} = O({\tau}^2),
\end{aligned}
\label{periodicN3}
\end{equation}
\begin{equation}
%\small
\begin{aligned}
\bar x = x^0 &+ \frac{{\tau}}{8}\left\{ 3f_1+f_3 + 4\alpha_2 (1-\alpha_2)(f_2-f_3)  \right\} \\
&+\frac{{\tau}^2}{48}\Bigl\{4L_{f_1}f_1+3 L_{f_1}f_3 +L_{f_3}f_3 + 6 \alpha_2 ( 2L_{f_1}f_2-2L_{f_1}f_3+L_{f_2}f_3-L_{f_3}f_3 ) \Bigr. \\ &+ 12\alpha_2^2(L_{f_1}f_3-L_{f_1}f_2 + L_{f_2}f_2-2L_{f_2}f_3 +L_{f_3}f_3)\Bigr\} + O({\tau}^3).
\end{aligned}
\label{barxN3}
\end{equation}
}}

{
{\bf Corollary~3.4.}
{\em Let the assumptions of Corollary~3.1 be satisfied with $N=4$, $\bar u_1=0$, $u_1^1=-u^2_1=-u^3_1=u^4_1$, and let $\tau_1,\tau_2,\tau_3,\tau_4>0$ be related to $\alpha_2, \alpha_3,\alpha_4$ by means of~\eqref{tau-alpha}.
Then
\begin{equation}
\tau_1=\left(\frac{1}{2}-\alpha_4\right)\tau,\; \tau_2 = \alpha_2\tau,\; \tau_3=\left(\frac{1}{2}-\alpha_2\right)\tau,\; \tau_4 = \alpha_4\tau,\;\; \alpha_2,\alpha_4\in \left(0,\frac{1}{2}\right),
\label{tauN4a}
\end{equation}
and
\begin{equation}
 \begin{aligned}
f_1&+f_3+\frac{\tau}{4}(L_{f_1}f_1-L_{f_3}f_3) + \frac{\tau^2}{24}(L_{f_1}^2 f_1 + L_{f_3}^2 f_3) \\
&+\alpha_2 \left\{2(f_2-f_3)+\tau(L_{f_1}f_2+L_{f_3}f_3) + \frac{\tau^2}{4}(L_{f_1}^2 f_2 -L_{f_3}^2 f_3  )\right\}\\
&-\alpha_4 \left\{2(f_1-f_4)+\tau(L_{f_1}f_1+L_{f_4}f_3) + \frac{\tau^2}{4} (L^2_{f_1}f_1 - L_{f_4}L_{f_3}f_3) \right\}\\
&+\alpha_2^2 \tau \left\{ L_{f_2}f_2 - L_{f_3}f_3 + \frac{\tau}{2}(L_{f_3}^2 f_3 + L_{f_1} L_{f_2} f_2) \right\} \\
&+ 2 \alpha_2 \alpha_4 \tau \left\{L_{f_4}f_3-L_{f_1}f_2 - \frac{\tau}{2}(L_{f_1}^2 f_2 + L_{f_4} L_{f_3} f_3)\right\}\\
&+\alpha_4^2 \tau \left\{ L_{f_1}f_1 - L_{f_4}f_4 + \frac{\tau}{2}(L_{f_1}^2 f_1 + L^2_{f_4} f_3)\right\} = O({\tau}^2),
\end{aligned}
\label{periodicN4}
\end{equation}
\begin{equation}
%\small
\begin{aligned}
\bar x &= x^0 + \frac{\tau}{2}\Bigl\{ \frac{1}{4}(f_1-f_3) +\alpha_2 (f_1+f_3)-\alpha_4 (f_1+f_4)+\alpha_2^2 (f_2-f_3)\Bigr. \\
& + \alpha_4 (\alpha_4-2\alpha_2) (f_1-f_4) \Bigl.\Bigr\}+ \frac{\tau^2}{4}\Bigl\{ \frac{1}{12}(L_{f_1}f_1+L_{f_3}f_3) + \frac{\alpha_2}{2}(L_{f_1}f_1-L_{f_3}f_3) \Bigr. \\
&- \frac{\alpha_4}{2}(L_{f_1}f_1-L_{f_4}f_3) + {\alpha_2^2}(L_{f_1}f_2 + L_{f_3}f_3)  + {\alpha_4^2}(L_{f_1}f_1 + L_{f_4}f_4) \\
& -2\alpha_2 \alpha_4 (L_{f_1}f_1 + L_{f_4}f_3) \Bigl.\Bigr\} +O({\tau}^3).
\end{aligned}
\label{barxN4}
\end{equation}
}}

{
Note that conditions~\eqref{tauN2},~\eqref{tauN3}, and~\eqref{tauN4a} follow from the assertion~\eqref{isotau-alpha} of Corollary~3.1 under the assumptions that $\bar u_1 =0$ and that the set $U$ is symmetric  in Corollaries~3.2--3.4.
Equations~\eqref{periodicN2},~\eqref{periodicN3}, and~\eqref{periodicN4} are obtained from the above Fliess expansions under the periodicity condition~\eqref{periodic_exp} with the $\tau_i$ given by~\eqref{tau-alpha}; and~\eqref{barxN2},~\eqref{barxN3},~\eqref{barxN4} are obtained by integrating the Fliess expansion for $x(t)$ over the period $[0,\tau]$.
}

\section{Numerical simulations and discussion}

The above analytical results will be applied in this section for computing switching controls in order to optimize the performance measure of the hydrolysis reaction {with the input reactant $\rm (CH_3CO)_2 O$  (denoted by $A$) and the product
$\rm CH_3 COOH$}.
Namely, we will treat system~\eqref{csys} with the vector fields given by~\eqref{f_notations} as a mathematical model of the chemical reaction
$\rm (CH_3CO)_2 O + H_2 O \to 2\, CH_3 COOH$ {with the following dimensionless parameters}:
\begin{equation}
{\bar n}=1,\;
{\varkappa=17.77,\; k_1 = 5.819 \cdot 10^7,\; k_2 = -8.99\cdot 10^5},\; \phi_1=\phi_2=1.
 \label{dimensionless_parameters}
\end{equation}
{
These values correspond to physical parameters of the adiabatic reaction by formulas from~\cite{NSP2015,CES2017}:}
{
$$
\begin{aligned}
& \varkappa = \frac{E_A}{R\bar T},\; R= 8.3144598\, \frac{J}{K\cdot mol},\; k_1 = k_0 {\bar C_A}^{\bar n-1} \frac{V}{\bar F},\; k_2 = \frac{\Delta H_R k_0 {\bar C_A}^{\bar n} V}{\rho c_p \bar T \bar F},\;  \\
& \phi_1 = \phi_2 = \frac{F}{\bar F},
\end{aligned}
$$
where $ E_A = 44.35\, \frac{kJ}{mol}$ is the activation energy, $k_0 = 1.4\cdot 10^5 \,s^{-1}$ is the collision factor,
$\Delta H_R = -55.5 \, \frac{kJ}{mol}$ is the reaction heat, $\rho c_p = 4.186\, \frac{kJ}{K \cdot l}$ is the product of the density and the heat capacity,
$V = 0.298\, l$ is the reactor volume, $F= 7.17\cdot 10^{-4}\,\frac{l}{s}$ is the volumetric flow-rate of the reaction stream, $\bar F =F$ is the steady-state flow-rate,
${\bar C_A}=0.3498\,\frac{mol}{l}$ is the steady-state outlet concentration of $A$, and $\bar T = 300.17\, K$ is the steady-state temperature in the reactor.
The state of system~\eqref{csys} is described by $x_1(t)=\frac{C_A(t)-\bar C_A}{\bar C_A}$ and $x_2(t)=\frac{T(t)-\bar T}{\bar T}$, where
the time variable $t$ corresponds to rescaling the physical time by $F/V$,
$C_A(t)$ is the concentration of $A$ in the reactor, and $T(t)$ is the temperature in the reactor.
}

{
The reaction is controlled by modulating the inlet concentration $C_{Ai}(t)$ of $A$ and the inlet temperature $T_i(t)$. These physical inputs correspond to the dimensionless controls $u_1(t)$ and $u_2(t)$ in~\eqref{csys}:
$$
\begin{aligned}
& u_1 (t)= \frac{k_1 (F-\bar F) e^{-\varkappa}}{\bar F} + \frac{(1+k_1 e^{-\varkappa})F}{\bar F \bar C_{Ai}}(C_{Ai}(t)-\bar C_{Ai}),\\
& u_2 (t)= \frac{k_2 (F-\bar F) e^{-\varkappa}}{\bar F} + \frac{(1+k_2 e^{-\varkappa})F}{\bar F \bar T_{i}}(T_{i}(t)-\bar T_{i}), \\
\end{aligned}
$$
where $\bar C_{Ai} = 0.74\, \frac{mol}{l} $ is the steady-state inlet concentration of $A$ and $\bar T_i = 295\, K,$ is the steady-state inlet temperature.
We assume the possibility of controlling the concentration $C_{Ai}(t)$ in the range of $(1\pm 0.85) \bar C_{Ai}$, and the temperature $T_i(t)$ in the range of $\bar T_i \pm 20\,K$.
This results in the control constraints $(u_1,u_2)^T\in U={\overline{\rm co}} \, U_b$ with
}
\begin{equation}
U_b = \left\{ \begin{pmatrix}u_1^{min} \\ u_2^{min} \end{pmatrix}, \begin{pmatrix}u_1^{min} \\ u_2^{max} \end{pmatrix},   \begin{pmatrix}u_1^{max} \\ u_2^{min} \end{pmatrix},  \begin{pmatrix}u_1^{max} \\ u_2^{max} \end{pmatrix}  \right\},
\label{strategy}
\end{equation}
and
{
\begin{equation}
u_1^{max} = -u_1^{min} = 1.798,\;  u_2^{max} = -u_2^{min} = 0.06663.
\label{maxsign}
\end{equation}
}

{
The equilibrium $x_1=x_2=0$ of system~\eqref{csys} with $u_1=u_2=0$ corresponds to the steady-state operating mode of the reactor with $C_A(t)=\bar C_A$ and $T(t)=\bar T$.
Our goal is to improve the conversion of $A$ to the product by using the same amount of input reactant $A$ over a period, which is formally stated as Problem~2.1 with the isoperimetric constant $\bar u_1 = 0$.
The application of condition~\eqref{D_cond} with parameters~\eqref{dimensionless_parameters} yields $D\approx 2.36 >0$.
Hence, we will follow the assumption stated in Remark~2.1 and consider bang-bang controls of the form~\eqref{u_tilde} with $N\le 4$ only.
}
%%%%%%%%%%%%%%%%%%%%%%%%%%%

\begin{figure}[h!]
{
\begin{minipage}[t]{.48\textwidth}
        \centering
        \includegraphics[width=1\linewidth,trim={2cm 14.8cm 8cm 2cm},clip]{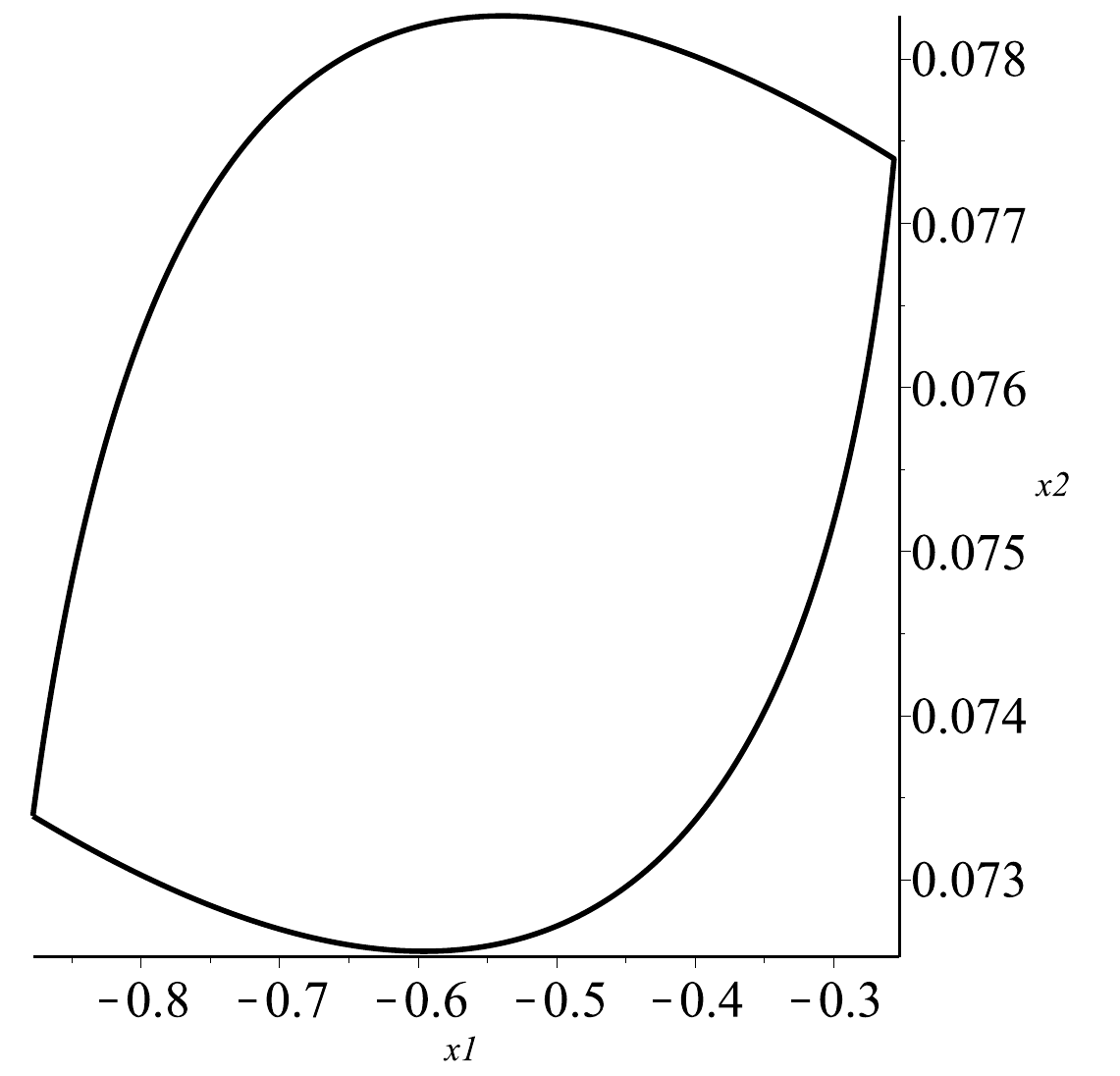}
        %\caption{Image 1.}%\label{fig:1a}
        {\scriptsize a: $N=2$, $J[x]\approx -0.566800$.}\\
        \vskip2ex
         \includegraphics[width=1\linewidth,trim={2cm 14.8cm 8cm 2cm},clip]{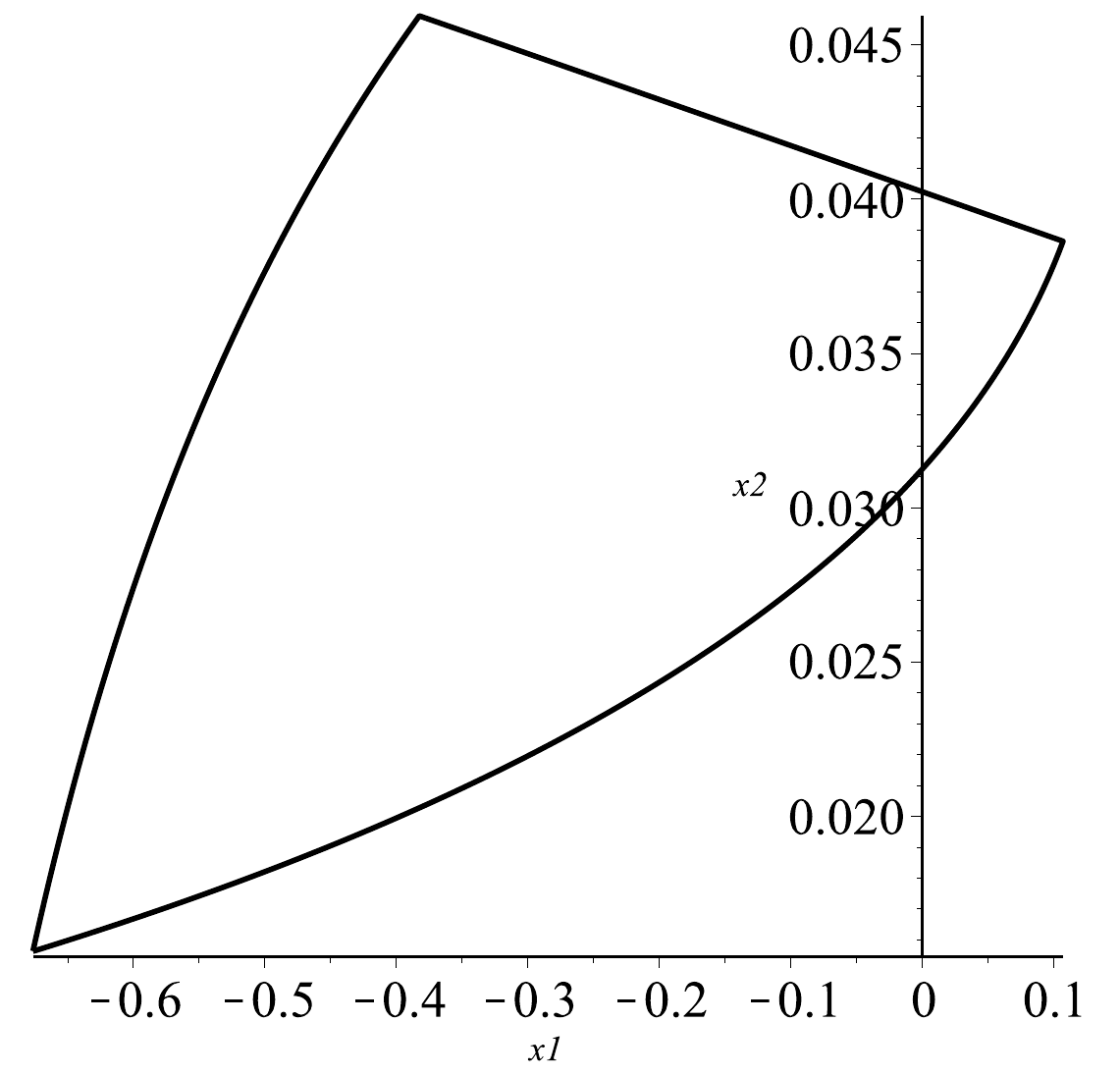}
        %\caption{Image 1.}%\label{fig:1a}
        {\scriptsize c:  $N=3$, $\alpha_2=0.2$, $J[x]\approx -0.287099$.}
    \end{minipage}
    \hfill
    \begin{minipage}[t]{.48\textwidth}
        \centering
        \includegraphics[width=1\linewidth,trim={2cm 14.8cm 8cm 2cm},clip]{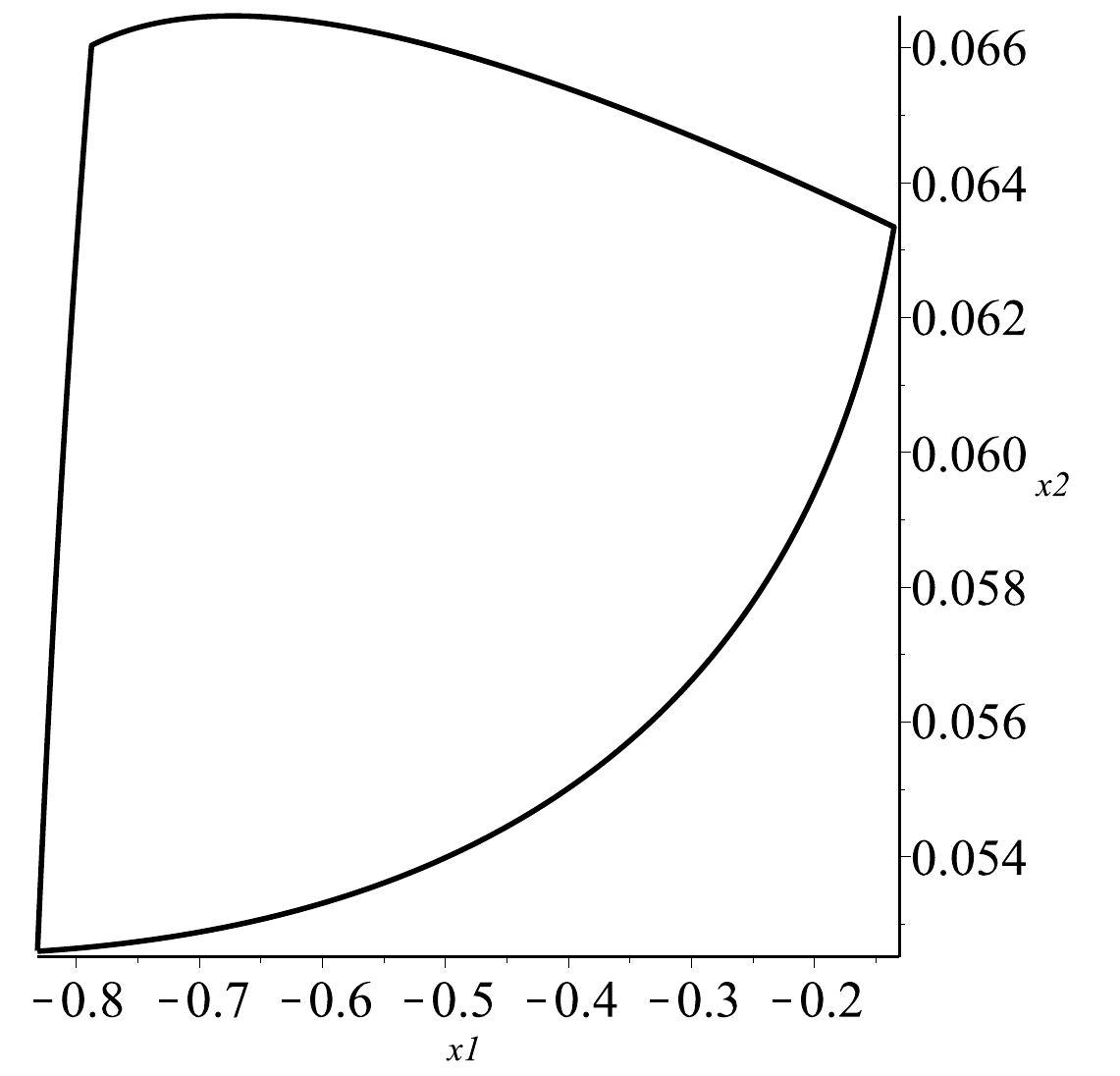}
        {\scriptsize b: $N=3$, $\alpha_2=0.4$, $J[x]\approx -0.482341$.}\\
        \vskip2ex
         \includegraphics[width=1\linewidth,trim={2cm 14.8cm 8cm 2cm},clip]{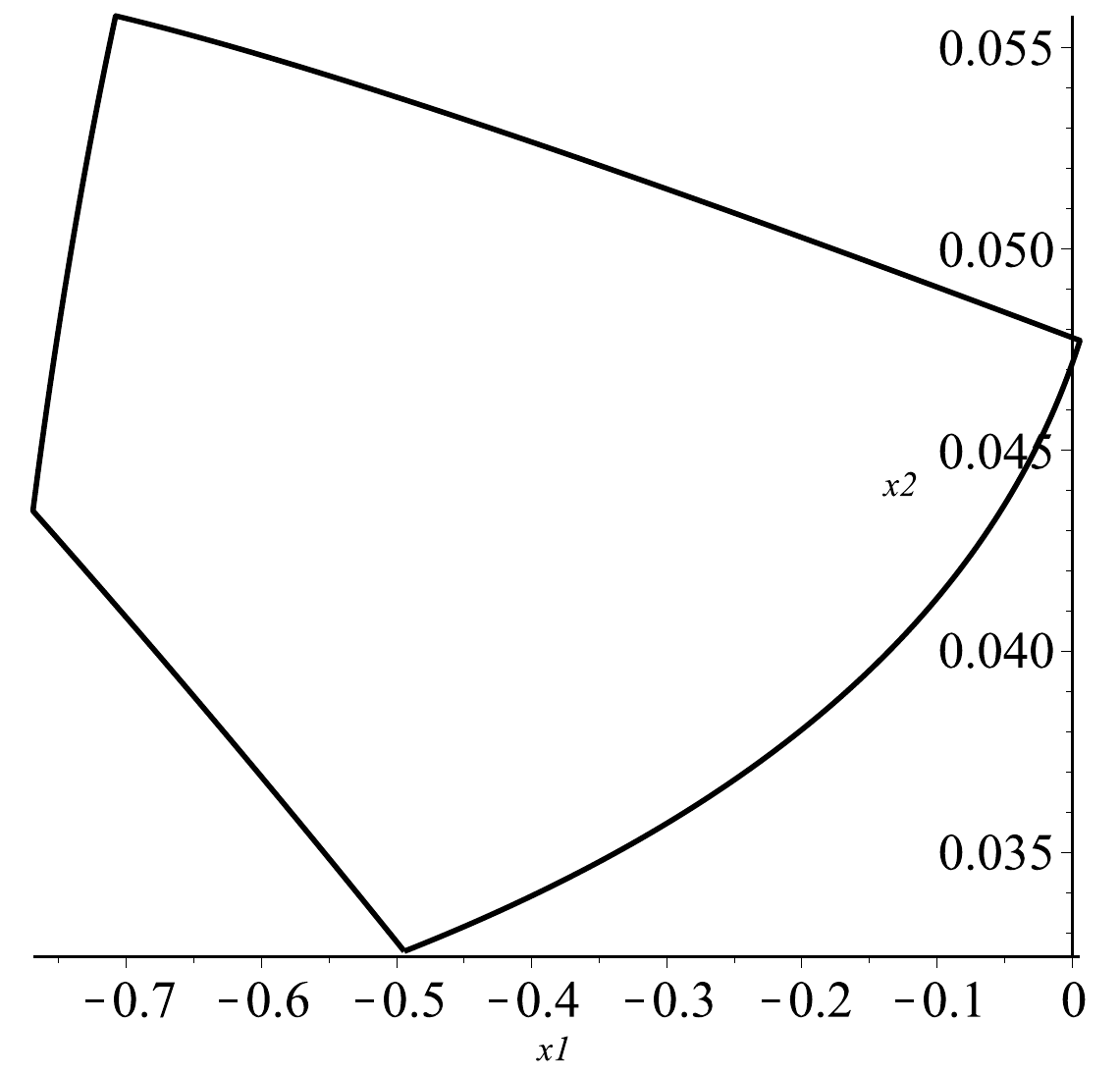}
        %\caption{Image 1.}%\label{fig:1a}
        {\scriptsize d: $N=4$, $\alpha_2=0.4$, $\alpha_4=0.1$, $J[x]\approx -0.379688$.}
        %\includegraphics[width=6cm,height=4cm,scale=1.6]{N2_t05_u1only.pdf}
        %\caption{Image 2.}%\label{fig:1b}
    \end{minipage}
\caption{Periodic trajectories of system~\eqref{csys} outside the equilibrium $x=0$, $\tau=1$.}
\label{fig:1}}
\end{figure}

{
Let us first consider constant controls. Note that system~\eqref{csys}  admits the following equilibria with $u_1=0$ and $u_2=\pm u_2^{max}$:
$$
\begin{aligned}
& x^- \approx (-0.566139,\, 0.075376 )^T\quad \text{with}\; u_1=0,\; u_2 = u_2^{max}, \\
& x^+ \approx (0.689896,\, -0.077288 )^T\quad \text{with}\; u_1=0,\; u_2 = -u_2^{max},
\end{aligned}
$$
and both of the above equilibria satisfy the constraints of Problem~2.1.	
The solution $x^-$ gives better performance in comparison with the trivial equilibrium:
\begin{equation}
 J[x^-]=x^-_1\approx -0.566139<0.
 \label{Jminus}
\end{equation}
%%%%
Hence, if a periodic trajectory $\{\tilde x(t)\}_{t\in [0,\tau]}$ of system~\eqref{csys} with some $\tilde u\in{\cal U}_{\tau}$ is contained in an $\varepsilon$-neighborhood of $x^-$ such that
 $\varepsilon\in (0,|x^-|)$, then
$$
J[\tilde x]\le x^-+\varepsilon < J[x]
$$
for any solution $x(t)\in X$ ($0\le t\le \tau$) of system~\eqref{csys} with $u\in {\cal U}_\tau$, where
$$
X = \{x\in\mathbb R^2 :\,x_1>x^-+\varepsilon\}.
$$
This means that the global solutions to Problem~2.1 cannot be obtained  by considering just small loops around $x=0$.
Some periodic trajectories of system~\eqref{csys} outside the equilibrium $x=0$ are shown in Fig.~\ref{fig:1} for controls~\eqref{u_tilde} with $N=2,3,4$.
These figures are obtained by numerical simulations in Maple.
}

{
From the practical viewpoint, the goal for studying Problem~2.1 is to optimize the performance in a neighborhood of the {\em given} steady-state $x=0$ by small variations of controls,
while the steady-state $x^-$ corresponds to {\em another} operating mode of the reactor (which {\em may not be desirable} due to requirements on the purity of the product or energy consumption).
To exclude the case $x=x^-$ from further consideration, we impose one more isoperimetric constraint:
\begin{equation}
\bar u_2= \frac{1}{\tau}\int_0^{\tau}u_2(t)dt=0.
\label{iso2}
\end{equation}
The above constraint corresponds to the assumption of using the same amount of energy as for the reference steady-state $x=0$.
Then the equilibrium $x^-$ is not a feasible solution anymore,  and we will study the admissible switching strategies for $N\le 4$.
}

Let the numbers ${\tau}>0$ and $N$ be fixed, and let the control $u (t)$ be given by~\eqref{u_tilde}:
\begin{equation}
u(t)=u^j \quad \text{for}\;\; t\in [t_{j-1},t_{j}),\quad j=\overline{1,N},
 \label{u_tilde2}
\end{equation}
for some switching scenario
$$
u^1,u^2,...,u^N \in U_b
$$
and switching times
\begin{equation}
0=t_0<t_1<...<t_N={\tau}.
\label{switching-t}
\end{equation}
{
It is easy to see that the cases $N=1$ and $N=3$ (if all $\tau_1$, $\tau_2$, $\tau_3$ are positive) are not consistent with the two isoperimetric constraints~\eqref{isotau} and~\eqref{iso2} as $\bar u_1=\bar u_2=0$ and $u_1^{max}>0$, $u_2^{max}>0$.
}

For $N\ge 2$, we express the switching times~\eqref{switching-t} using~\eqref{tau-alpha} in terms of the positive parameters
$\alpha_2$, ...., $\alpha_N$ such that $\alpha_2+...+\alpha_N<1$:
\begin{equation}
\begin{aligned}
t_1& = \left(1-\sum_{j=2}^N \alpha_j\right){\tau},\\
t_j& = t_{j-1} + \alpha_j {\tau},\quad j=2,..., N.
\end{aligned}
\label{tau-alpha2}
\end{equation}

%%%%%%%%%%%%%%%%%%%%%%%%%%%%%%%%%%%%%%%%%%%%

%By assuming that the comment of Subsection~\ref{subsection21} concerning the maximal number of switchings remains valid for the nonlinear system~\eqref{csys} in a neighborhood of the origin,
%we may consider controls of the form~\eqref{u_tilde} with $N\le 5$ only.

Let us now consider the control~\eqref{u_tilde2} with $N=2$:
\begin{equation}
u(t) = \left\{\begin{matrix}u^1, & t\in [0,(1-\alpha_2) {\tau}), \\ u^2, & t\in [(1-\alpha_2) {\tau},{\tau}]. \end{matrix}\right.
\label{control2}
\end{equation}
Then $u(t)$ satisfy the isoperimetric constraint~\eqref{isotau} {with $\bar u_1=0$ and~\eqref{iso2}} if and only if
\begin{equation}
{u^2 = -u^1},\; \alpha_2 = 1/2,
\label{condN2}
\end{equation}
as it follows from {Corollary~3.2}.
Note that, for an arbitrary initial condition $x(0)=x^0\in {\mathbb R}^2$, the solution $x(t)$ of system~\eqref{csys} with~\eqref{control2} is not necessary ${\tau}$-periodic.
Thus we apply Corollary~{3.2} to find a relation between $x^0$ and ${\tau}$  such that the corresponding solution $x(t)$ of~\eqref{csys} with the chosen switching strategy satisfies the periodic boundary
condition $x(0)=x({\tau})$.
{
We rewrite formulas~\eqref{periodicN2} and~\eqref{barxN2} in terms of the original vector fields $f_0$, $g_1$, $g_2$ of system~\eqref{csys},~\eqref{f_notations} as follows:
\begin{equation}
f_0+ \frac{{\tau}}{4}L_{\tilde g_1}f_0+ \frac{{\tau^2}}{24}\left(L_{f_0}L_{f_0}+L_{\tilde g_1}L_{\tilde g_1}\right)f_0 = O(\tau^2),
\label{periodicN2g}
\end{equation}
\begin{equation}
\bar x = x^0+\frac{\tau}{4}\tilde g_1 + \frac{\tau^2}{24}L_{f_0}f_0 + O(\tau^3),
\label{xbarN2}
\end{equation}
where we have assumed that $\tilde g_1 = u^1_1 g_1 + u^1_2 g_2$ is a constant vector field.
Note that the vector fields in~\eqref{periodicN2g} and~\eqref{xbarN2} are evaluated at $x=x^0$, so that the equation~\eqref{periodicN2g} implicitly defines the map
\begin{equation}
x^0 = c_1 \tau + c_2 \tau^2 + O(\tau^3)
\label{x0tau}
\end{equation}
with
$$
c_1 = - (G_{x^0})^{-1} G_\tau,\; c_2 = -\frac{1}{2} (G_{x^0})^{-1}  \left\{ c_1^T G_{x^0 x^0} c_1 + G_{\tau x^0} c_1 + G_{\tau \tau} \right\},
$$
for small $\tau>0$, where $G(x^0,\tau)$ denotes the left-hand side of~\eqref{periodicN2g}, and
$G_\tau$, $G_{x^0}$, $G_{\tau x^0}$, $G_{x^0 x^0}$ are corresponding derivatives of $G$ at $\tau=0$, $x^0=0$ (we treat $G_{x^0}$ as the Jacobian matrix and $G_{x^0 x^0}$ as the Hessian).
We have:
$$
G_{x^0} = \begin{pmatrix}-\phi_1 - k_1 e^{-\varkappa} & - k_1 \varkappa e^{-\varkappa} \\ - k_2 e^{-\varkappa}  & -\phi_2 - k_2 \varkappa e^{-\varkappa} \end{pmatrix} \approx \begin{pmatrix} -2.115 & -19.820 \\ 0.0172 & -0.694\end{pmatrix}.
$$
We see that the matrix $G_{x^0}$ is nonsingular, ${\rm det}(G_{x^0})\approx 1.809\neq 0$, so that the proposed method for approximate computation of the inital conditions for periodic trajectories with $N=2$ and small periods succeeds:
$c_1 \approx (-0.4495,-0.0167)^T$ and $c_2 \approx (-0.14133,-0.00218)^T$ in formula~\eqref{x0tau}.
}

%%%%%%%%%%%%%%%%%%%%%%%%%%%%%%%%%%%%%%%%%%%%%%%%%%%%%%%%%%%%%%%%%%%%%%%%%%%%%%%%%%%%%%%%%%%%%
\begin{figure}[h!]
{
\begin{minipage}[t]{.31\textwidth}
        \centering
        \includegraphics[width=1\linewidth,trim={2cm 19cm 12cm 2cm},clip]{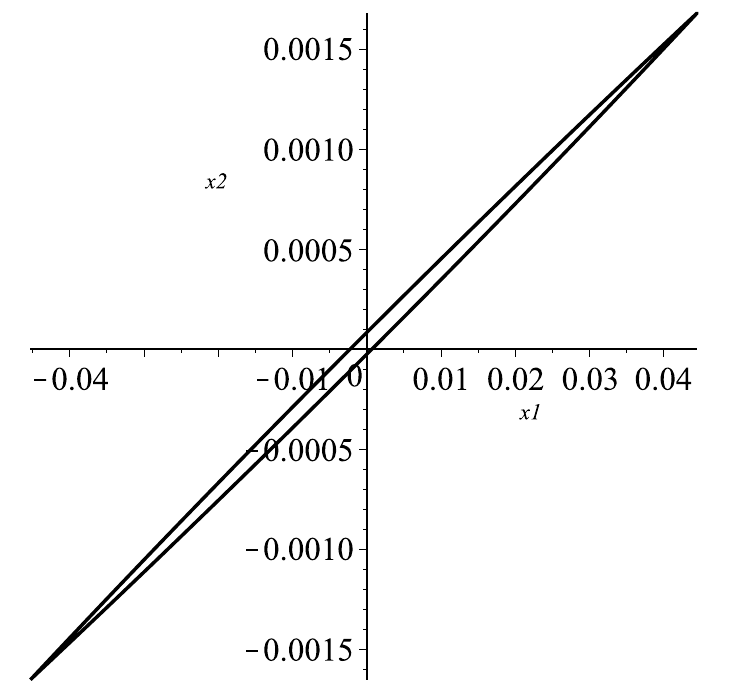}
        {\scriptsize a: $\tau=0.1$, $J[x]\approx -0.00040$.}\\
        \vskip2ex
         \includegraphics[width=1\linewidth,trim={2cm 14.8cm 8cm 2cm},clip]{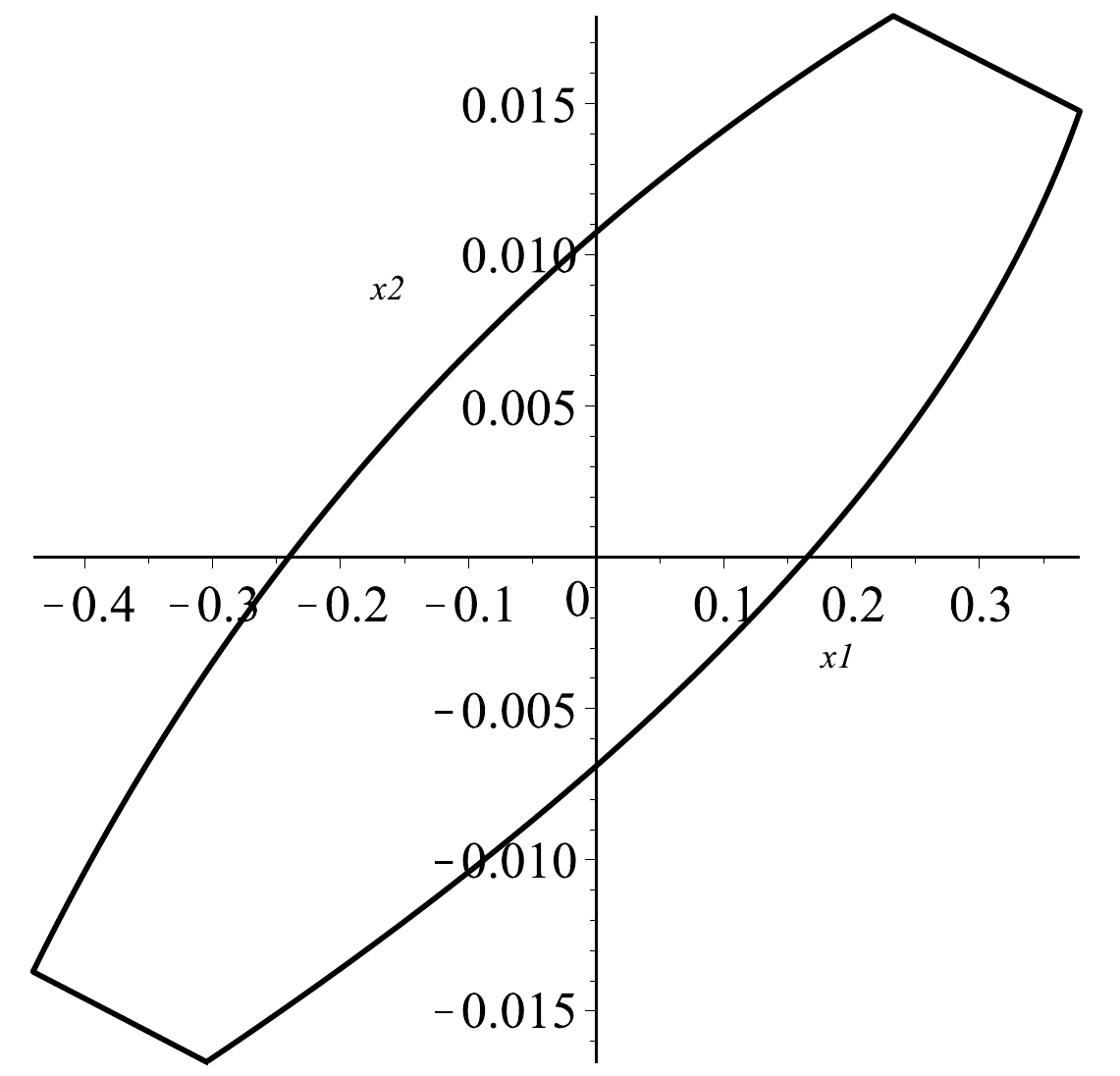}
        {\scriptsize d: $\tau=1$, $\alpha_2=\alpha_4=0.05$, $J[x]\approx -0.03112$.}
    \end{minipage}
    \hfill
    \begin{minipage}[t]{.31\textwidth}
        \centering
        \includegraphics[width=1\linewidth,trim={2cm 19cm 12cm 2cm},clip]{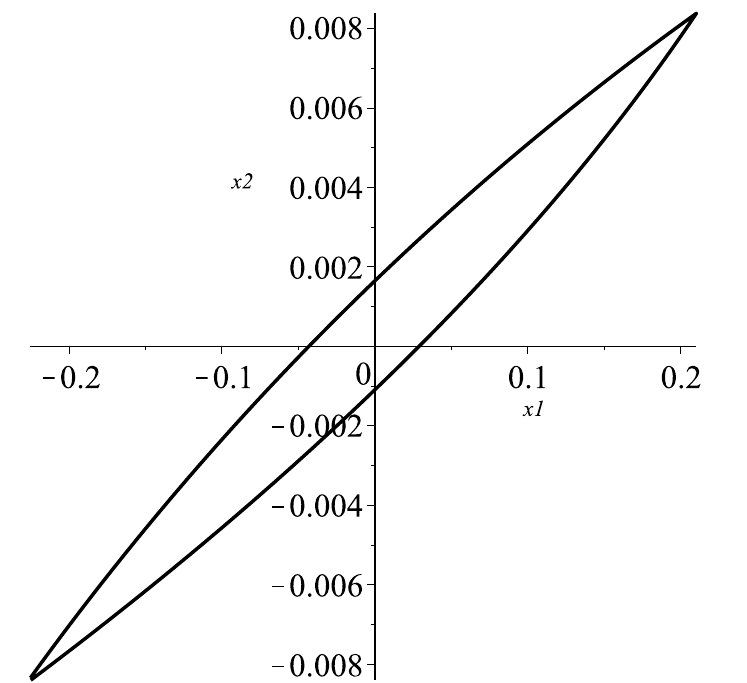}
        {\scriptsize b: $\tau=0.5$, $J[x]\approx  -0.00726$.}\\
                \vskip2ex
         \includegraphics[width=1\linewidth,trim={2cm 14.8cm 8cm 2cm},clip]{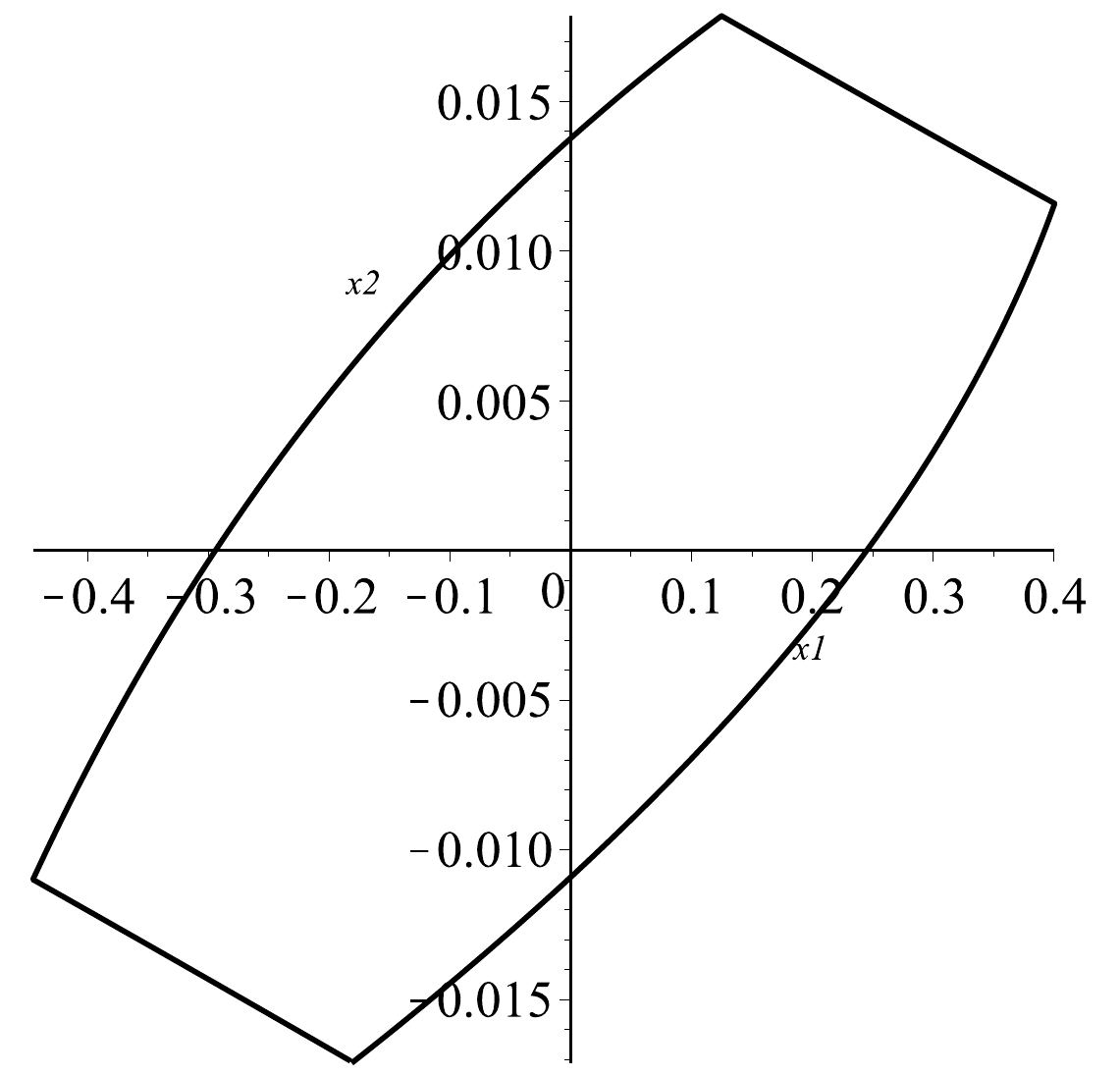}
        {\scriptsize e: $\tau=1$, $\alpha_2=\alpha_4=0.1$, $J[x]\approx -0.02497$.}
    \end{minipage}
        \hfill
    \begin{minipage}[t]{.31\textwidth}
        \centering
        \includegraphics[width=1\linewidth,trim={2cm 19cm 12cm 2cm},clip]{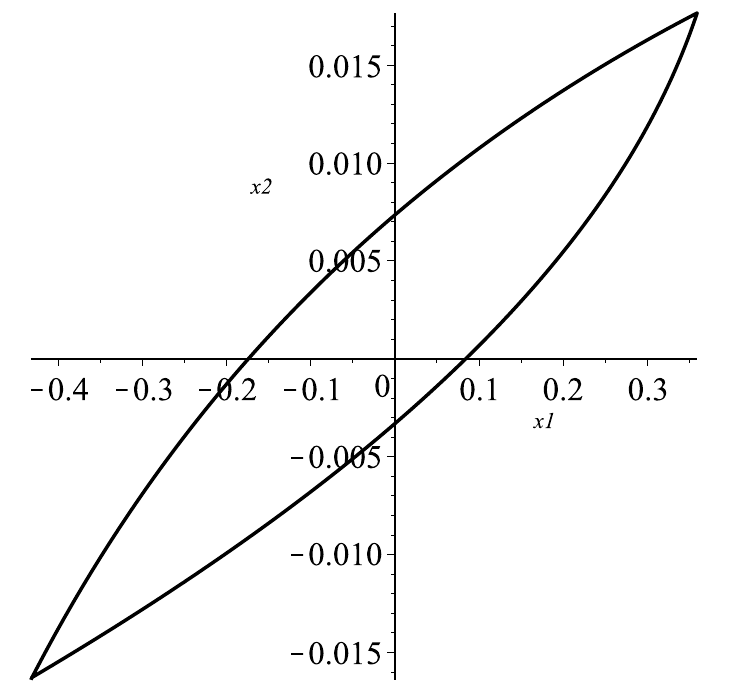}
        {\scriptsize c: $\tau=1$, $J[x]\approx -0.03385$.}\\
        \vskip2ex
         \includegraphics[width=1\linewidth,trim={2cm 14.8cm 8cm 2cm},clip]{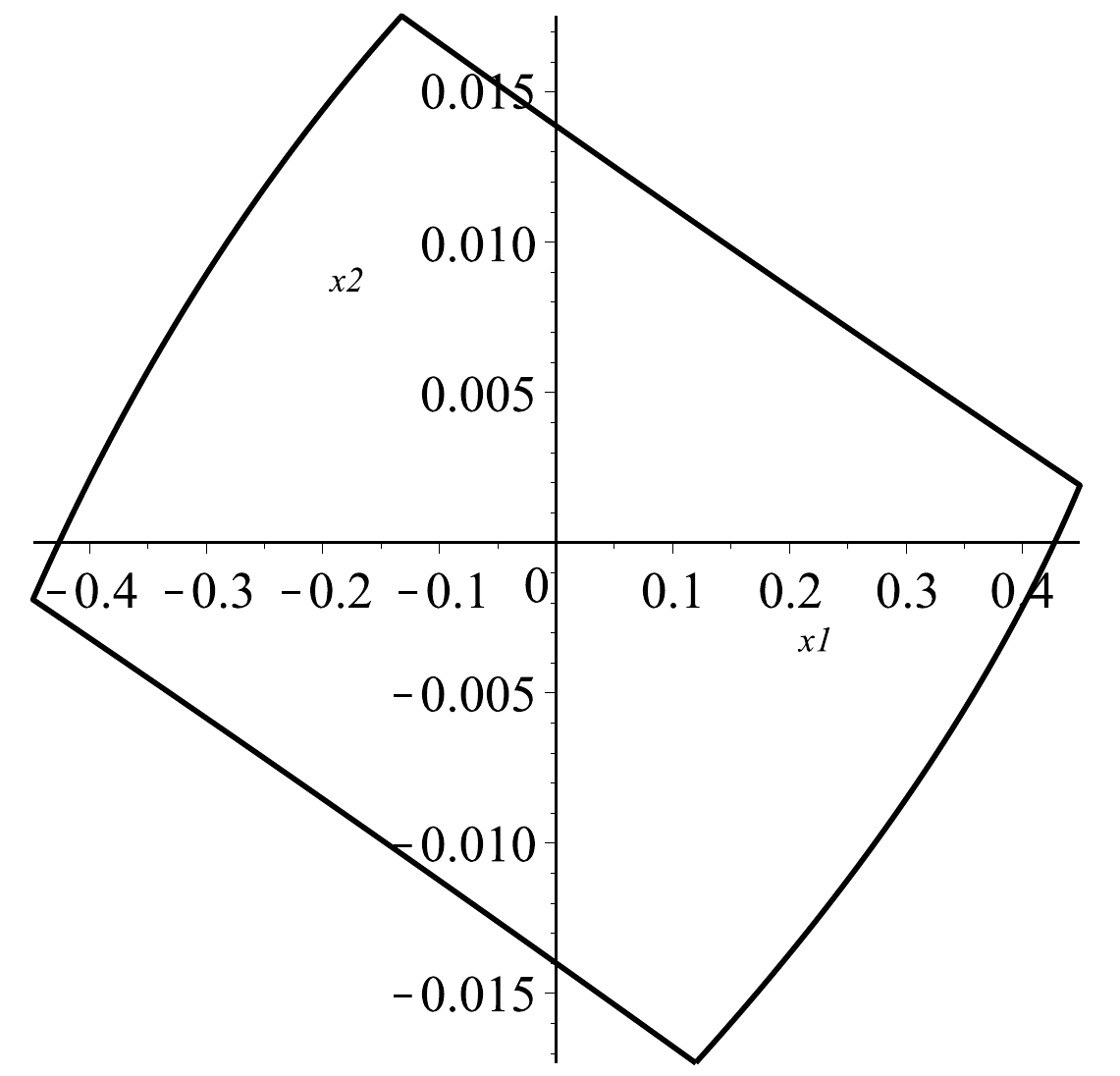}
        {\scriptsize f: $\tau=1$, $\alpha_2=\alpha_4=0.25$, $J[x]\approx -0.00295$.}
    \end{minipage}
\caption{Periodic trajectories around $x=0$ with small $\tau$: $N=2$ (a--c) and $N=4$ (d--f).}
\label{fig:2}}
\end{figure}

Trajectories of system~\eqref{csys} with the control~\eqref{control2} {under the switching strategy $u^1 = -u^2 = (u_1^{max}, u_2^{max})^T$ are shown in Fig.~\ref{fig:2} (a--c)}.
%These illustrations show that the corresponding solutions $x(t)$ are ${\tau}$-periodic.
{We also summarize the obtained numerical results in Table~1 for different values of $\tau$.}
 {As one can see,  $J[x]<0$ along the constructed periodic trajectories}, so that the  control~\eqref{control2} improves the reactor performance in comparison to the steady-state trajectory $x=0$.

  \begin{table}[h!]
\begin{center}
{
\begin{tabular}{|c|c|c|c|c|}
\hline
$\tau$ & ${x^0}$ & $J[x]$ & $J^{est}_2(x^0,\tau)$\\ \hline
0.1 & $(-0.04529, -0.00165)^T$ & -0.00040 & -0.00039 \\ \hline
0.2 & $(-0.09130, -0.00325)^T$ & -0.00188 & -0.00186 \\ \hline
0.3 & $(-0.13601, -0.00498)^T$ & -0.00276 & -0.00269 \\ \hline
0.4 & $(-0.18170, -0.00659)^T$ & -0.00554 & -0.04935\\ \hline
0.5 & $(-0.22511, -0.00840)^T$ & -0.00726 & -0.00730 \\ \hline
0.6 & $(-0.27032, -0.00986)^T$ & -0.01248 & -0.01268 \\ \hline
0.7 & $(-0.31291, -0.01152)^T$ & -0.01674 & -0.01728 \\ \hline
0.8 & $(-0.35420, -0.01315)^T$ & -0.02172 & -0.02269 \\ \hline
0.9 & $(-0.39371, -0.01482)^T$ & -0.02709 & -0.02868 \\ \hline
1.0 & $(-0.43202, -0.01630)^T$ & -0.03385 & -0.03580 \\ \hline
\end{tabular}
\caption{Simulation results for $N=2$.}}
\end{center}
\label{tab:01}
\end{table}

The significance of the above theoretical results is also underpinned by possible applications of the formula~\eqref{xbarN2} for analytic approximation of the cost $J$.
Indeed, {let us denote by $J^{est}_2(x^0,\tau)$ the first coordinate of $\bar x$ in~\eqref{xbarN2}
with higher order terms being neglected.
%without the terms of order $O(\tau^3)$.
The values of $J^{est}_2(x^0,\tau)$ are presented in Table~1 together with the corresponding cost, and we see that $J^{est}_2(x^0,\tau)$ gives a good approximation of $J[x]$ for $\tau\le 1$.
Moreover, by substituting~\eqref{x0tau} into $J^{est}_2(x^0,\tau)$ and computing its Taylor expansion at $\tau=0$, we obtain
\begin{equation}
J^{est}_2 = c^*\cdot \tau^2 + O(\tau^3),
 \label{J2est}
\end{equation}
where $c^* \approx -0.141$ for the considered example. As $c^*<0$, we conclude that {\em any} periodic trajectory corresponding to the control~\eqref{control2} is profitable in comparison to $x=0$, provided that $\tau>0$ is small enough.
}

To {study} the behavior of system~\eqref{csys} with  $N=4$, we choose the control~\eqref{u_tilde2} with the
 following switching scenario:
 $$
 u^1 = -u^3 = (u_1^{max}, u_2^{max})^T, \; u^2 = -u^4 =(u_1^{min}, u_2^{max})^T.
 $$
 Then {Corollary~3.4 implies, assuming that $\tilde g_1 = u^1_1 g_1 + u^1_2 g_2 = {\rm const}$ and $\tilde g_2 = u^2_1 g_1 + u^2_2 g_2 = {\rm const}$:
 \begin{equation}
  \begin{aligned}
   f_0 &+ (\alpha_2-\alpha_4)(\tilde g_1+\tilde g_2) + \frac{\tau}{2}\Bigl\{ \frac{1}{2}L_{\tilde g_1}f_0 + 2\alpha_2 L_{f_0}f_0 -\alpha_4 (2 L_{f_0} + L_{\tilde g_1}  -L_{\tilde g_2} )f_0\Bigr.\\
   & + (\alpha_2-\alpha_4)^2(L_{\tilde g_1}+L_{\tilde g_2})f_0\Bigl.\Bigr\} + \tau^2\Bigl\{ \frac{1}{24}(L^2_{f_0} + L^2_{\tilde g_1})\Bigr. \\
   & + \frac{\alpha_2}{4}( L_{f_0}L_{\tilde g_1} + L_{\tilde g_1} L_{f_0} ) - \frac{\alpha_4}{8}( 2 L_{f_0} L_{\tilde g_1} + L_{\tilde g_1}L_{f_0} +  L_{\tilde g_2}L_{f_0} + L_{\tilde g_1}^2 +L_{\tilde g_2}L_{\tilde g_1} )\\
   &+\frac{\alpha_2^2}{4}(2 L_{f_0}^2 - L_{f_0} L_{\tilde g_1}+ L_{f_0} L_{\tilde g_2} + L_{\tilde g_1}^2 + L_{\tilde g_1} L_{\tilde g_2})\\
   & -\frac{\alpha_2 \alpha_4}{2} (2 L_{f_0}^2 + L_{\tilde g_1}L_{f_0} - L_{\tilde g_2}L_{f_0} + L_{\tilde g_1}^2 + L_{\tilde g_2} L_{\tilde g_1}) \\
   &+ \frac{\alpha_4^2}{2}(2 L_{f_0}^2 + L_{f_0} L_{\tilde g_1} -L_{f_0} L_{\tilde g_2} + L_{\tilde g_1} L_{f_0} - L_{\tilde g_2} L_{f_0}+ L_{\tilde g_1}^2 + L_{\tilde g_2}^2)\\
   & + \frac{\alpha_2^3}{6}(L_{f_0} L_{\tilde g_1}+L_{f_0} L_{\tilde g_2}+L_{\tilde g_1}L_{f_0} +L_{\tilde g_2}L_{f_0}-L_{\tilde g_1}^2+L_{\tilde g_2}^2)\\
   &- \frac{\alpha_2^2\alpha_4}{2}(L_{f_0} L_{\tilde g_1}+L_{f_0} L_{\tilde g_2}+L_{\tilde g_1}L_{f_0}+L_{\tilde g_2}L_{f_0})\\
   &+\frac{\alpha_2\alpha_4^2}{2}(L_{f_0} L_{\tilde g_1}+L_{f_0} L_{\tilde g_2}+ L_{\tilde g_1}L_{f_0}+L_{\tilde g_2}L_{f_0} - L_{\tilde g_2}^2 )\\
   &-\frac{\alpha_4^3}{6} (L_{f_0} L_{\tilde g_1}+ L_{f_0} L_{\tilde g_2}+L_{\tilde g_1}L_{f_0}+ L_{\tilde g_2}L_{f_0}+ L_{\tilde g_1}^2- L_{\tilde g_2}^2)\Bigl. \Bigr\}f_0=O(\tau^3),
   \label{periodicN4g}
 \end{aligned}
 \end{equation}
 \begin{equation}
 \begin{aligned}
  \bar x &= x^0+ \frac{\tau}{2}\Bigl\{ \frac12 \tilde g_1 + 2\alpha_2 f_0 -\alpha_4 (2f_0 + \tilde g_1 - \tilde g_2)+(\alpha_2 - \alpha_4)^2(\tilde g_1 + \tilde g_2)\Bigr\} \\
  &+ \frac{\tau^2}{2}\Bigl\{ \frac{1}{12} L_{f_0} f_0 + \frac{\alpha_2}{2}L_{\tilde g_1} f_0 - \frac{\alpha_4}{4}(L_{\tilde g_1} + L_{\tilde g_2})f_0 + \alpha_2^2   L_{f_0}f_0 \Bigr. \\
  & + \alpha_4 (\alpha_4-2\alpha_2) (L_{f_0}+\frac{1}{2} L_{\tilde g_1} - \frac{1}{2}L_{\tilde g_2})f_0 \Bigl.\Bigr\} + O(\tau^3).
 \end{aligned}
\label{xbarN4}
\end{equation}
 }

{
Similarly to the previous consideration~\eqref{x0tau}, we exploit the periodicity condition~\eqref{periodicN4g} to define $x^0$ for small values of $\tau$:
\begin{equation}
x^0 = \bar c_1(\alpha_2,\alpha_4) \tau + \bar c_2 (\alpha_2,\alpha_4)\tau^2 + O(\tau^3),
\label{x0tau4}
\end{equation}
where $\bar c_i(\alpha_2,\alpha_4)$ are polynomials of $(\alpha_2,\alpha_4)$.
The above $x^0$ is used for computing periodic trajectories of system~\eqref{csys}, depending on the parameters $\alpha_2,\alpha_4\in (0,1/2)$ and $\tau>0$.
Fig.~\ref{fig:2} (d--f) illustrate these trajectories for controls of the form~\eqref{u_tilde2} with $N=4$.
We observe that $J[x]<0$ for the considered trajectories; thus, the proposed controls with $N=4$ improve the performance in comparison to the steady-state operation {$x=0$}.
}

{
To compare the reactor performances with $N=2$ and $N=4$, we denote by $J^{est}_4(x^0,\tau,\alpha_2,\alpha_4)$ the first coordinate of $\bar x$ in~\eqref{xbarN4} without the terms $O(\tau^3)$
and substitute the expression~\eqref{x0tau4} into $J^{est}_4(x^0,\tau,\alpha_2,\alpha_4)$. As a result, we have
\begin{equation}
J^{est}_4 = \bar c(\alpha_2,\alpha_4)\cdot \tau^2 + O(\tau^3),
\label{Jest4}
\end{equation}
where $\bar c(\alpha_2,\alpha_4)$ is a polynomial of $(\alpha_2,\alpha_4)$. Our numerical study shows that $\bar c(\alpha_2,\alpha_4)>c^*$ for all $(\alpha_2,\alpha_4)\in (0,1/2)^2$, where the constant $c^*$ appears in~\eqref{J2est}.
}

{
Thus, for the considered numerical example, we have analyzed the asymptotic expansions $J_2^{est}$ and $J_4^{est}$ under the periodicity conditions~\eqref{periodicN2g} and~\eqref{periodicN4g} for small $\tau$
to conclude that the controls with $N=4$ do not improve the performance in comparison with the case $N=2$ locally,
along the  periodic trajectories near zero. This conclusion is also confirmed by the simulations presented in Fig.~\ref{fig:2}.
}

\begin{figure}[hb!]
{
\begin{minipage}[t]{.31\textwidth}
        \centering
        \includegraphics[width=1\linewidth,trim={2cm 19cm 12cm 2cm},clip]{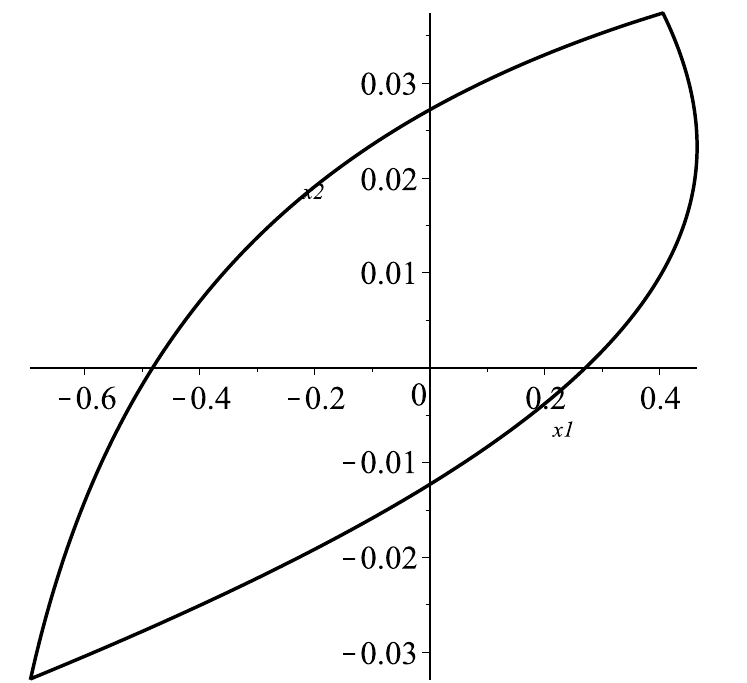}
        {\scriptsize a: $\tau=2$, $J[x]\approx -0.10898$.}\\
        \vskip2ex
         \includegraphics[width=1\linewidth,trim={2cm 14.8cm 8cm 2cm},clip]{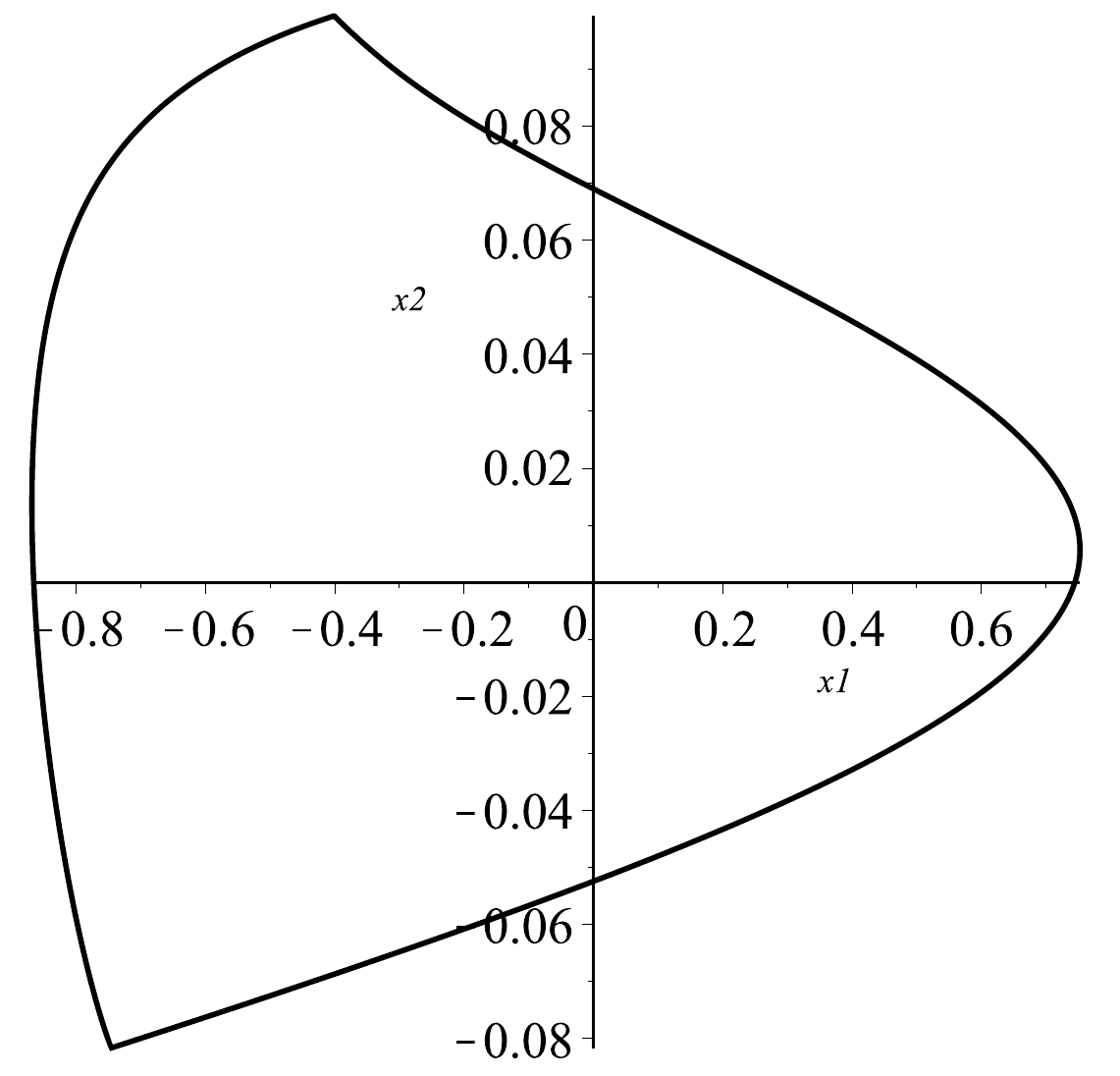}
        {\scriptsize d: $\tau=10$, $J[x]\approx -0.41555$.}
    \end{minipage}
    \hfill
    \begin{minipage}[t]{.31\textwidth}
        \centering
        \includegraphics[width=1\linewidth,trim={2cm 14.8cm 8cm 2cm},clip]{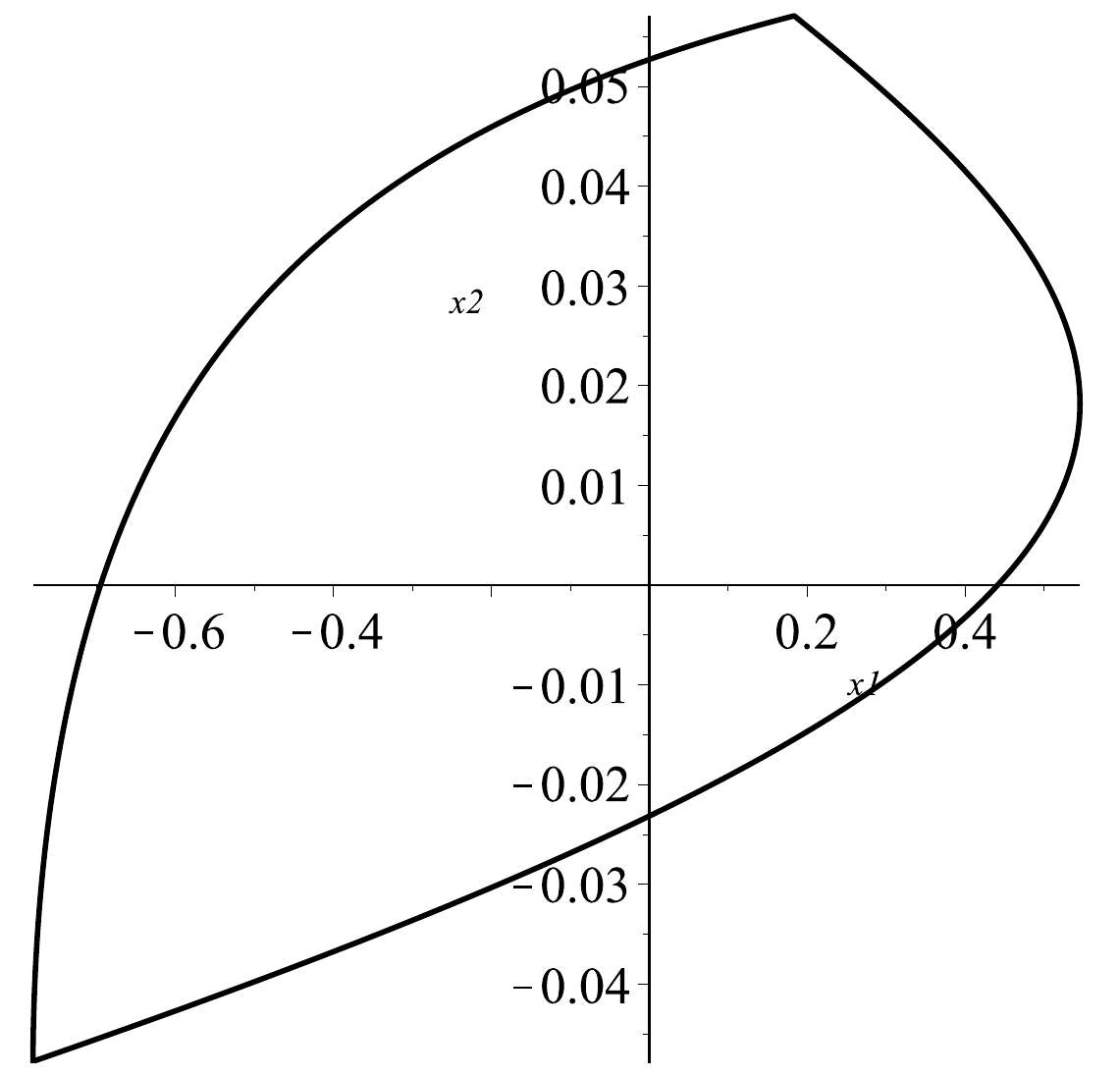}
        {\scriptsize b: $\tau=3$, $J[x]\approx -0.18622$.}\\
                \vskip2ex
         \includegraphics[width=1\linewidth,trim={2cm 14.8cm 8cm 2cm},clip]{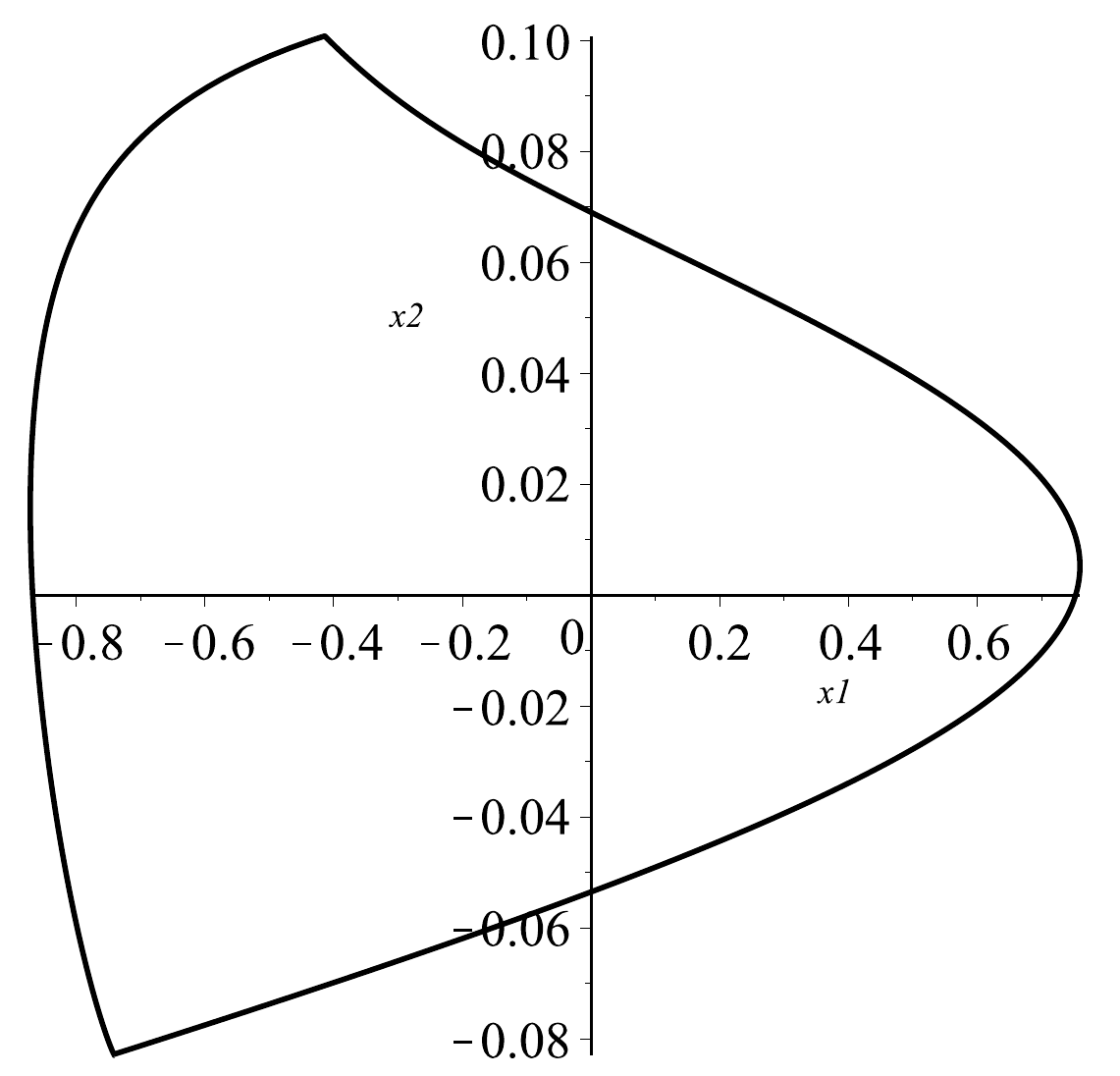}
        {\scriptsize e: $\tau=100$, $J[x]\approx -0.56096$.}
    \end{minipage}
        \hfill
    \begin{minipage}[t]{.31\textwidth}
        \centering
        \includegraphics[width=1\linewidth,trim={2cm 14.8cm 8cm 2cm},clip]{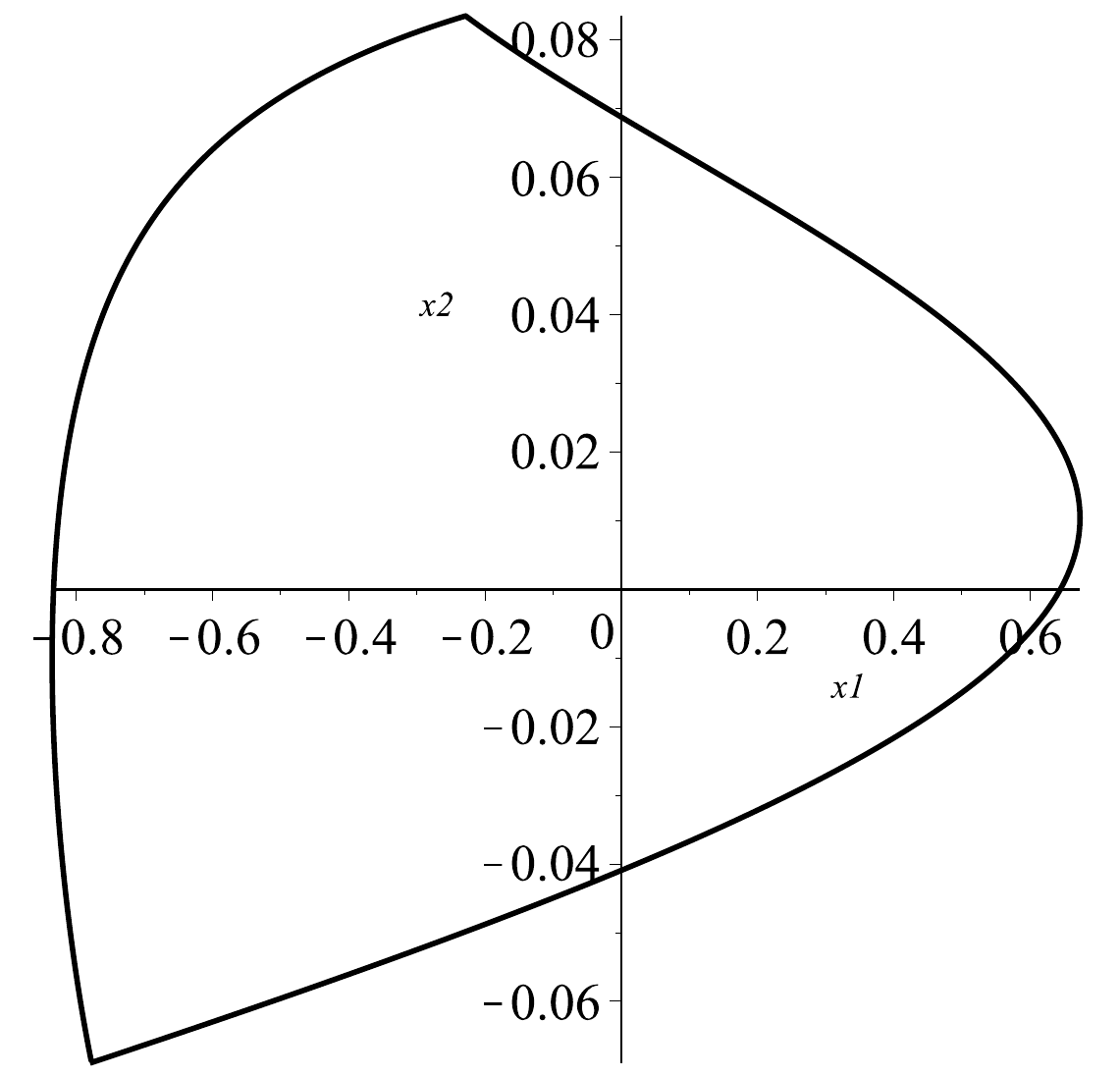}
        {\scriptsize c: $\tau=5$, $J[x]\approx -0.28761$.}\\
        \vskip2ex
         \includegraphics[width=1\linewidth,trim={2cm 14.8cm 8cm 2cm},clip]{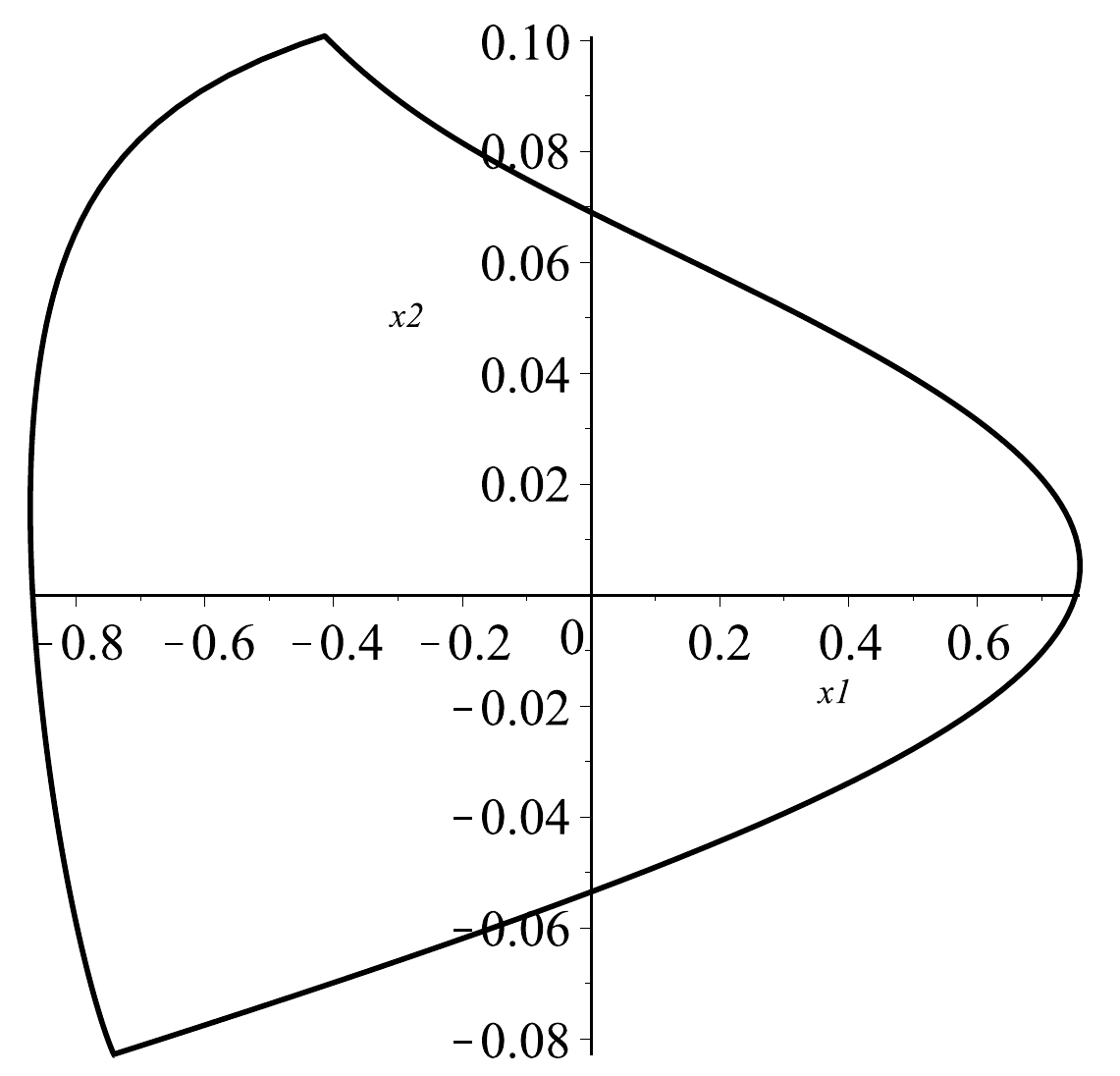}
        {\scriptsize f: $\tau=1000$, $J[x]\approx -0.57565$.}
    \end{minipage}
\caption{Periodic trajectories around $x=0$ with large $\tau$, $N=2$.}
\label{fig:x}}
\end{figure}

{
The above analytical results are valid for small values of $\tau$.
However, the data of Table~1 and estimates of the type~\eqref{J2est},~\eqref{Jest4} suggest that the increasing of $\tau$ leads to the decreasing of the cost under appropriate control strategies.
To illustrate this behavior, we compute periodic trajectories of system~\eqref{csys} numerically with the control~\eqref{control2} for increasing values of $\tau$.
These trajectories together with their costs are presented in Fig.~\ref{fig:x}.}

%%%%%%%%%%%%%%%%%%%%%%%%%%%%%%%%%%%%%%%%%%%%%%%%%%%%%%%%%%%%%%%%%%%%%%%%%%%%%%%%%%%%%%%%%%%%%

{
We observe that the cost $J[x]$ is monotonically decreasing when $\tau$ is increasing, so that the sequence of controls~\eqref{control2} with $\tau\to+\infty$ may be considered as a candidate for a minimizing sequence for Problem~2.1.
The corresponding trajectories $\{x(t)\}_{t\in [0,\tau]}$ converge very fast (in the orbital sense) to some limit curve, as  the plots for $\tau=10$, $\tau=100$, and $\tau=1000$ look almost identical in Fig.~\ref{fig:x}.
We also note that the trajectory with $\tau=1000$  ensures even better performance $J[x]\approx -0.57565$ than the steady-state solution $x=x^-$ with $J[x^-]\approx -0.566139$.
However, such trajectories with a large time horizon exhibit considerable deviations from the reference steady-state $x=0$, which may not be acceptable in practical applications.
}

\section{Conclusions}
The proposed control design scheme generalizes the results of~\cite{DSTA2017} for multidimensional nonlinear control-affine systems and bang-bang strategies with an arbitrary number of switchings.
As it follows from the comparison of analytical and numerical results in Section~4,
our approach can be used for estimating the performance of nonlinear chemical reactors analytically and
improving the conversion ``$A\to$ product'' with respect to the trivial steady-state solution.
It should also be noted that this approach has a potential for minimizing the cost $J$ analytically by varying the {phase} parameters $\alpha_j$ {under a fixed time horizon $\tau$}
in order to tune the phases of different input signals in an optimal fashion.
{The development of analytical tools for estimating performance measures with large values of $\tau$ remains to be an issue for future study.}

\section*{Acknowledgements}
This work was supported by the Strategic Innovation Fund of the Max Planck Society.




%\bibliographystyle{elsarticle-num}
%\bibliography{switching}


%\begin{thebibliography}{00}


%\bibitem{}

%\end{thebibliography}

\end{document}